\documentclass[12pt,a4paper, draft]{amsart}
\usepackage[english]{babel}
\usepackage{amssymb}
\usepackage{amscd}
\usepackage{mathrsfs}
\usepackage{amssymb,amsmath,amsthm}
\usepackage{mathtools}

\usepackage[all]{xy}
\usepackage{tikz-cd}
\usepackage{pgfplots,tikz}
\usetikzlibrary{decorations.pathmorphing,patterns}
\usepackage{afterpage}
\usepackage{hyperref}
\usepackage{verbatim}
\usepackage{url}

\usepackage{geometry}
\newtheorem{thm}{Theorem}[section]
\newtheorem{cor}[thm]{Corollary}
\newtheorem{lem}[thm]{Lemma}
\newtheorem{prop}[thm]{Proposition}
\newtheorem{ex}[thm]{Example}

\newtheorem{prob}{Problem}

\newcommand{\RK}[1]{\mathcal{R}^{#1}\mathfrak{K}}
\newcommand{\RH}[1]{\mathcal{R}^{#1}\operatorname{Hom}}

\theoremstyle{remark}

\theoremstyle{definition}

\newcommand{\fK}{\mathfrak K}

\newcommand{\co}{\operatorname{co}}

\newcommand\R{\mathbb R}
\newcommand\C{\mathbb C}
\newcommand\K{\mathbb K}

\newcommand{\CAT}{\textbf}
\newcommand{\To}{\longrightarrow}
\newcommand{\OP}{\mathfrak}

\renewcommand\Im{\operatorname{Im}}
\DeclareMathOperator{\coker}{coker}
\newcommand{\Ext}{\operatorname{Ext}}
\newcommand{\Hom}{\operatorname{Hom}}

\newcommand{\PO}{\operatorname{PO}}
\newcommand{\PB}{\operatorname{PB}}

\newcommand\restr[2]{{
		\left.\kern-\nulldelimiterspace 
		#1 
		\right|_{#2} 
}}
\newcommand{\norm}[1]{{\left\vert\kern-0.25ex\left\vert #1
		\right\vert\kern-0.25ex\right\vert}}
\newcommand{\nnorm}[1]{{\left\vert\kern-0.25ex\left\vert\kern-0.25ex\left\vert #1
		\right\vert\kern-0.25ex\right\vert\kern-0.25ex\right\vert}}

\title[Compact derivation in Banach spaces]{Derived functors associated to the ideal of compact operators in Banach spaces}

\thanks{The research of the first three authors has been supported in part by project PID2023-146505NB-C21 funded by MCIU/AEI/10.13039/5011000110331/FEDER/UE. The research of the first two authors has been supported in part by projects IB24002 and GR24056 funded by Junta de Extremadura. The third author has also been supported by project PID2024-162214NB-I00 funded by MCIN/AEI/10.13039/501100011033The fourth author benefited from a predoctoral grant associated to ``Programa de Financiaci\'on UCM-Banco Santander (CT24/25)'', awarded by Universidad Complutense de Madrid.}

\author[Cabello]{F\'elix Cabello S\'anchez}
\address{Instituto de Matem\'aticas\\ Universidad de Extremadura\\ Avenida de Elvas\\ 06071-Badajoz\\ Spain}
\email{fcabello@unex.es}

\author[Castillo]{Jes\'us \,M.\,F. Castillo}
\address{Instituto de Matem\'aticas\\ Universidad de Extremadura\\ Avenida de Elvas\\ 06071-Badajoz\\ Spain}
\email{castillo@unex.es}

\author[Salguero]{\\ Alberto Salguero-Alarc\'on}
\address{Universidad Complutense de Madrid\\ Plaza de las Ciencias, s/n\\ 28040-Madrid\\ Spain}
\email{albsalgu@ucm.es}

\author[Trejo]{Nazaret Trejo-Arroyo}
\address{Universidad Complutense de Madrid\\ Plaza de las Ciencias, s/n\\ 28040-Madrid\\ Spain}
\email{ ntrejo@ucm.es}

\numberwithin{equation}{section}


\begin{document}
	\maketitle

\begin{abstract} We compute the derived functors of (the functors associated to) the ideal of compact operators in Banach spaces and obtain new results about the extension and lifting of compact operators.\end{abstract}

\section{Introduction}

\subsection{Purpose}
The extension and lifting of operators, especially that of compact operators, in Banach spaces has been a topic present since the very beginning of the theory. The first impulse came from Fredholm theory and integral equations, but that quickly developed, at least since Lindenstrauss' memoir \cite{memoir}, through its connections with the local theory of Banach spaces.

In this paper we undertake the study of the ideal of compact operators from a homological perspective focusing on the bifunctor $(X,Y) \rightsquigarrow \frak K(X,Y)$, its associated functors $X \rightsquigarrow
\frak K(X,Y)$ and $Y \rightsquigarrow
\frak K(X,Y)$ and their derived functors.

Our research is addressed to Banach spacers and homological people alike. So we find the Gordian knot of how making the paper readable for both groups since Banach techniques might probably not be in the toolkit of categorical people, and at the same time the subject makes unavoidable a certain amount of ``abstract nonsense" that Banach spacers might find hard to digest. To circumvent this issue
we organized the paper making every effort to explain that the abstract part of this paper is not just ``turning the crank", as Lefschetz might have been afraid of \cite{lef}, warning in advance that Banach spaces do not form an Abelian category precisely because topology plays an essential role.

\subsection{Summary} The paper is organized as follows. This section contains the general introduction and explains some notational conventions used throughout the paper. Section 2 displays some background material on compact operators, Yoneda $\Ext$ for Banach spaces and derived functors, mainly to fix the notation. Section~\ref{sec:RnK} presents the derived functors of $X \rightsquigarrow \frak K(X,Y)$ (using projectives) and $Y \rightsquigarrow \frak K(X,Y)$ (using injectives).
The ideal of compact operators is so well balanced that
  the $n$-th derived of $\frak K(-,Y)$ at $X$ agrees with the $n$-th derived of $\frak K(X,-)$ at $Y$ and, in the end, we have just one (graded, bi-) functor that, quite naturally, we denote by $\Ext_\fK^n$. Like all derived functors, $\Ext_\fK^n$ comes with two long homology sequences, which we describe in detail.

In Section~4 we move to a more Banach-oriented level to study the (obvious) natural transformation $\eta: \Ext_{\fK}^n\to \Ext^n$.
We isolate and characterize \emph{compact $n$-extensions}, that is, $n$-extensions in the image of $\eta$
as those that can be obtained from a projective  resolution (equivalently, from an injective resolution) using compact operators. Focusing on the first derived functors, we show, by means of examples and results proved ad hoc, that the canonical map $\Ext_\fK(X,Y)\to\Ext(X,Y)$ is neither injective or surjective. In spite of this fact there are large families of Banach spaces for which the vanishing of  $\Ext_\fK(X,Y)$ implies the vanishing of $\Ext(X,Y)$, or viceversa. In a final turn of screw we summon quasilinear maps to provide an intrinsic, topological characterization of compact extensions and show that (under rather mild hypotheses) they all live
in the non-Hausdorff part of $\Ext$.

The closing Section~\ref{sec:misc}
presents a number of problems related to this research  that we we find interesting and, believe it or not, we have
not been able to solve.

\subsection{Notational conventions}

We consider linear spaces over $\R$ or $\C$; and $\mathbb K$ denotes the ground field.  A homogeneous map $B:X\to Y$, acting between normed or quasinormed spaces, is said to be bounded if there is a constant $C$ such that $\|B(x)\|\leq C\|x\|$ for every $x\in X$. The map is said to be compact if the image of the unit ball of $X$ is a relatively compact set of $Y$. Operators are always assumed to be linear and bounded, otherwise we speak of linear maps. An embedding is an operator that is an isomorphism onto its range. Given Banach spaces $X, Y$, we denote by $\mathfrak L(X,Y)$ the space of linear operators from $X$ to $Y$, while $\mathfrak K(X,Y)$ denotes the subspace of compact operators and $\mathfrak F(X,Y)$ the finite rank operators. If $Y=X$ we just write $\frak L(X)$, and the same for $\fK$ and $\frak F$. Compact operators are often denoted by lowercase greek letters: $\alpha, \beta, \tau, \sigma\dots $ Throughout the paper ${\bf B}$ is the category of Banach spaces and operators and ${\bf V}$ is the category of vector spaces and linear maps. These, and the underlying category of sets and maps, are the only ones we will need. Given a map $f:A\to B$ acting  between sets (possibly with additional structures) we use $f^*$ (respectively, $f_*$) to indicate composition with $f$ on the right (respectively, on the left) in a variety of situations. The identity map on $A$ is denoted by ${\mathbf 1}_A$.

\section{Background}

\subsection{Exact sequences of Banach spaces} A finite or infinite sequence of Banach spaces and operators
$\xymatrix{ \dots \ar[r] & U_{k-1} \ar[r]^{u_{k-1}} & U_k  \ar[r]^-{u_{k}} & U_{k+1} \ar[r] & \dots
  }$ (typically denoted by $\mathscr U$) is exact at $U_k$ if the kernel of $u_k$ coincides with the range of $u_{k-1}$. A sequence is said to be exact if it is exact at every position.   If one merely has $u_k u_{k-1}=0$ for all $k$, then $\mathscr U$ is called a complex. A morphism of complexes $\mathscr U \to \mathscr V$ is a collection of arrows making a commutative diagram
$$
\xymatrixcolsep{4pc}
\xymatrix{(\mathscr U) &
\cdots \ar[r] & U_{k-1} \ar[r]^{u_{k-1}} \ar[d] & U_k  \ar[r]^{u_{k}}  \ar[d]& U_{k+1} \ar[r] \ar[d] & \dots
 & \\
(\mathscr V) &
\cdots \ar[r] & V_{k-1} \ar[r]^{v_{k-1}} & V_k  \ar[r]^{v_{k}} & V_{k+1} \ar[r] & \dots& }$$
A covariant (resp. contravariant) additive functor $F: {\bf B}\to {\bf V}$ is said to be left-exact if for every short exact sequence
$0\to A \to B \to C \to 0$, the sequence
$0\to FA \to FB \to FC$ (resp. $0\to FC \to FB \to FA$) is exact. A functor that preserves short exact sequences is said to be exact.

\subsection{A reminder of Yoneda Ext for Banach spaces}
Let $X$ and $Y$ be Banach spaces. An \emph{$n$-extension} of $X$ by $Y$ is an exact sequence  in {\bf B} of the form
	\begin{equation}\tag{$\mathscr Z$}\label{dia:YZnZ1X}
		\begin{tikzcd}
			0 \arrow[r] & Y \arrow[r] & Z_n \arrow[r] & \cdots \arrow[r] & Z_1 \arrow[r] & X \arrow[r] & 0
		\end{tikzcd}
	\end{equation}
Within the class of all $n$-extensions of $X$ by $Y$ we consider the equivalence relation generated by  the order that declares $\mathscr Z \gtrdot \mathscr Z' $ (also written $\mathscr Z' \lessdot \mathscr Z $)  if and only if there is a commutative diagram
\begin{equation}\label{dia:equi in Ext n}
\xymatrixrowsep{1.5pc}
\xymatrix{
& & Z_n \ar[r] \ar[dd] & \cdots \ar[r] & Z_1 \ar[dd] \ar[dr] &\\
0\ar[r] &Y \ar[ur] \ar[dr] &&&& X\ar[r] &0\\
& & Z_n' \ar[r] & \cdots \ar[r] & Z_1' \ar[ur]
}
\end{equation}
It can be shown \cite[Corollary 6.40]{frer-sieg} that $\mathscr Z$ is equivalent to $\mathscr Z'$ (written $\mathscr Z\sim \mathscr Z'$) if and only if there exists $\mathscr Z''$
such that either $\mathscr Z \gtrdot \mathscr Z'' \lessdot  \mathscr Z'$ or  $\mathscr Z \lessdot \mathscr Z'' \gtrdot  \mathscr Z'$. We denote by $\Ext^n(X,Y)$ the set of equivalence classes of $n$-extensions of $X$ by $Y$. If $n=1$ we will often omit it. Extensions (that is, $1$-extensions) are just short exact sequences $\begin{tikzcd}[cramped, column sep=1.25em]
		0 \arrow[r] & Y \arrow[r] & Z \arrow[r] & X \arrow[r] & 0
	\end{tikzcd}$, whose meaning is that $Y$ is isomorphic to a closed subspace of $Z$ with $Z/Y$ isomorphic to $X$. In this case equivalence is  simpler because if $n=1$ the (only) vertical arrow in \eqref{dia:equi in Ext n} must be an isomorphism, by the three lemma and the open mapping theorem.\medskip	
	
Given an $n$-extension $\mathscr Z$ and an operator  $u: X' \to X$, the associated \emph{pullback $n$-extension} $\mathscr Z u$ is the lower row of the diagram
	\begin{equation}
		\notag
		\begin{tikzcd}
	&&		0 \arrow[r] & Y \arrow[r] \arrow[d, equal] & Z_n \arrow[r] \arrow[d, equal] & \cdots \arrow[r] & Z_1 \arrow[r, "q"]                            & X \arrow[r]                 & 0 & (\mathscr{Z}) \\
	&&		0 \arrow[r] & Y \arrow[r]                                & Z_n \arrow[r]                                                                & \cdots \arrow[r] & \PB \arrow[r] \arrow[u] & X' \arrow[r] \arrow[u, "u"] & 0 & (\mathscr Z u)
		\end{tikzcd}
	\end{equation}
	where $\PB=\lbrace (z,x') \in Z_1 \times X': q(z)= u(x') \rbrace$, and the arrows $\PB \to Z_1$ and $\PB \to X'$ are the restrictions of the canonical projections. Dually, the \emph{pushout $n$-extension} $v \mathscr Z $ associated to an operator $v: Y \to Y'$ is the lower row of the diagram
	\begin{equation}
		\notag
		\begin{tikzcd}
		&&	0 \arrow[r] & Y \arrow[r, "j"] \arrow[d, "v"] & Z_n \arrow[r] \arrow[d] & \cdots \arrow[r] & Z_1 \arrow[r]                                & X \arrow[r]                                & 0 & (\mathscr Z) \\
		& &	0 \arrow[r] & Y' \arrow[r]             & \PO \arrow[r]                                                   & \cdots \arrow[r] & Z_1 \arrow[r] \arrow[u, equal] & X \arrow[r] \arrow[u, equal] & 0 & (v \mathscr Z)
		\end{tikzcd}
	\end{equation}
Here $\PO$ is the quotient of the direct sum $Y' \oplus Z_n$ by the closure of the subspace $\lbrace (v(y), -j(y)) : y \in Y \rbrace$, and the arrows $Y' \to \PO$ and $Z_n \to \PO$ are the compositions of the natural injections into $Y' \oplus Z_n$ with the quotient map. Pullbacks and pushouts provide the assignment $(X, Y) \rightsquigarrow\Ext^n(X,Y)$ with a structure of set-valued bifunctor, contravariant in $X$, covariant in $Y$. Actually $\Ext^n(X,Y)$ can be given a natural linear structure (the Baer sum) that makes $\Ext^n$ a ${\bf V}$-valued functor. For $n=1$ the zero of $\Ext^n(X,Y)$ is (the class of) the direct sum extension $0\to Y\to Y\oplus X\to X\to 0$, and for $n>1$
is the class of the $n$-extension
	\begin{equation}
		\notag
		\begin{tikzcd}
			0 \arrow[r] & Y \arrow[r, equal] & Y \arrow[r, "0"] & \cdots \arrow[r, "0"] & Z \arrow[r, equal] & Z \arrow[r] & 0
		\end{tikzcd}
	\end{equation}

	\subsection{Compact operators} An operator $\tau:X\to Y$, acting between Banach spaces, is said to be \emph{compact} if it takes the unit ball of $X$ into a relatively compact subset of $Y$. Following Pietsch \cite{pietsch}, we will denote $\mathfrak K$ the operator ideal of compact operators. The ideal $\mathfrak K$ is \emph{injective}, with the meaning that whenever $\imath$ is an embedding and $\imath \tau\in \mathfrak K$ then $\tau\in \mathfrak K$; and it is also \emph{surjective}, with the meaning that whenever $\pi$ is a quotient map then $\tau\pi \in \mathfrak K$ implies $\tau \in \mathfrak K$. Given a Banach space $E$, let us consider the contravariant functor  $\mathfrak K(-,E):  \CAT{B}  \to  \CAT{V}$, which assigns to each Banach space $X$ the vector space $\mathfrak K(X,E)$ of all compact operators from $X$ to $E$ and to each operator $u:X \to Y$ the map $u^*:\mathfrak K(Y,E) \to \mathfrak K(X,E)$ given by $u^*(\tau) = \tau u$. Analogously, the covariant functor $\mathfrak K(E,-)$ assigns to each Banach space $X$ the vector space $\mathfrak K(E,X)$ and to each operator $u:X \to Y$ the map $u_*: \mathfrak K(E,X) \to \mathfrak K(E,Y)$ given by $u_*(\tau) = u \tau$. Both functors $\mathfrak K(E,-)$ and $\mathfrak K(-,E)$ are left-exact. Folklore facts from Banach space theory (see \cite[Theorem 4.2]{lpspaces} and \cite[Theorems 5 and 6]{mediterranean}) can be reformulated as:
	\begin{prop}\label{prop:K-is-exact} Given a Banach space $E$:
		\begin{enumerate}
			\item $\mathfrak K(-,E)$ is exact if and only if $E$ is an $\mathscr{L}_{\infty}$-space.
			\item $\mathfrak K(E, -)$ is exact if and only if $E$ is an $\mathscr{L}_1$-space.
		\end{enumerate}
	\end{prop}

\subsection{Derivation of functors in Banach spaces} Before going any further, let us clarify the context in which we deal with derived functors. Let us begin observing  {\bf B} is not an Abelian category: an operator may be both monic (injective) and epic (with dense range) without being an isomorphism. This causes some major and some minor complications that can be overcome
by completely adhering to the easily accessible notes by Frerick and Sieg \cite{frer-sieg} and considering the category of Banach spaces as an \emph{exact category} in the sense of \cite[Definition 3.5]{frer-sieg} with its maximal exact structure \cite[Chapter~4]{frer-sieg}. This just means that the admissible exact sequences are, well, all exact sequences. See also \cite{buhler:ex, buhler:alg}.

\medskip

Thus, as it is customary in Banach space theory:
\begin{itemize}
\item A Banach space $I$ is injective if $\mathfrak L(-,I)$ is exact.
\item A Banach space $P$ is projective if $\mathfrak L(P,-)$ is exact.
\end{itemize}
That is, $I$ is injective if, for every Banach space $Z$ and every subspace $Y\subset Z$, every operator from $Y$ to $I$ admits an extension to $Z$; and $P$ is projective if for every quotient map $Z\to X$ every operator from $P$ to $X$ has a lifting to $Z$. Every projective Banach space is isomorphic to $\ell_1(\Gamma)$ for some index set $\Gamma$ (a classical result by K\"othe \cite{kothe}, see also \cite{trimming}). On the other hand, while it is clear that the spaces $\ell_\infty(\Gamma)$ are injective, a complete description of injective spaces is unknown.

Nevertheless, every injective Banach space $I$ is an $\mathscr L_\infty$-space and so $\mathfrak K(-,I)$ is exact. Dually, if $P$ is projective, then $\mathfrak K(P, -)$ is exact because $P$
is of type $\mathscr L_1$.

The category of Banach spaces has \emph{enough} injectives (every Banach space is isomorphic to a closed subspace of some injective space)
and \emph{enough} projectives (every Banach space is a quotient of some projective space). This implies that every Banach space $Y$ admits an \emph{injective resolution}: an exact sequence
\begin{equation}\label{eq:inj-res-X}
\xymatrix{
0 \ar[r] & Y \ar[r] & I_1 \ar[r] & I_2 \ar[r] & \dots}
\end{equation}
in which each $I_k$ is injective and any Banach space $X$ has a \emph{projective resolution}: an exact sequence
\begin{equation}\label{eq:pro-res-X}
\xymatrix{
\dots \ar[r] & P_2 \ar[r] & P_1 \ar[r] & X  \ar[r] & 0}
\end{equation}
in which each $P_k$ is projective.\medskip

Assume that $F$ is a covariant left-exact functor from {\bf B} to an Abelian category, that for every purpose in this paper will be {\bf V}. Then the right-derived functors of $F$ arise as the homology of the complex that one obtains by applying $F$ to an injective resolution (\ref{eq:inj-res-X}) of $Y$, namely
$$
(\mathcal R^n F)Y = \dfrac{\ker(FI_{n+1}\to FI_{n+2})}{\Im(FI_{n}\to FI_{n+1})}
	$$
Since $F$ is left-exact, it transforms kernels into kernels, so $(\mathcal R^0 F)Y= \ker(FI_{1}\to FI_{2}) = FY$. When $F$ is a contravariant left-exact functor the right-derived functors arise as the homology of the complex one obtains by applying $F$ to a projective resolution (\ref{eq:pro-res-X}) of $X$; namely
$$
(\mathcal R^n F)X = \dfrac{\ker(FP_{n+1}\to FP_{n})}{\Im(FP_{n+2}\to FP_{n+1})}
	$$

Applying this recipe to the covariant functor $\OP L(X, -)$ or to the contravariant functor $\OP L(-,Y)$ associated to the bifunctor $(X,Y) \rightsquigarrow \frak L(X,Y)$, it turns out that the derived functors of $\OP L(-, Y)$ at $X$ using projectives agree with the derived functors of $\OP L(X, -)$ at $Y$ using injectives and with $\Ext^n(X,Y)$. Namely,
\begin{equation}\label{eq:RnL=Extn}
\Ext^n(X,Y) = \mathcal{R}^n\mathfrak L (X,-)\,Y  = \mathcal{R}^n\mathfrak L (-,Y)\, X .
\end{equation}

Even if a complete proof can be seen in \cite[Remark 6.44]{frer-sieg}. Since the understanding of these identities is of paramount importance for our approach, we will briefly explain them in detail right now.\medskip
	
Let us begin with the natural transformation of contravariant functors $\mathcal{R}^n\mathfrak L (-,Y) \to \Ext^n(-,Y)$. Fix a projective resolution of $X$ and set $K_n=\ker(P_n\to P_{n-1})$. The sequence $0 \to K_n \to P_n \to \dots \to P_1 \to X \to 0$ is exact and each $v \in \frak L(K_n, Y)$ induces the pushout sequence
$$\xymatrix{0 \ar[r] & K_n \ar[r] \ar[d]_v & P_n \ar[r] \ar[d] & \cdots \ar[r] & P_1 \ar[r]& X \ar[r]& 0 \\
0 \ar[r] & Y \ar[r]& \PO \ar[r] & \cdots \ar[r] & P_1 \ar[r] \ar@{=}[u]& X \ar[r] \ar@{=}[u] & 0 }
$$
The identity in (\ref{eq:RnL=Extn}) at the object level encodes the following facts:
\begin{itemize}
	\item $v$ extends to $P_n$ if and only if the lower sequence is zero in $\Ext^n(X,Y)$.
	\item Every $n$-extension of $X$ by $Y$ is equivalent to one that arises in  this way for some operator $K_n\to Y$.
\end{itemize}	
	
Dually, the other identity in (\ref{eq:RnL=Extn}) means, among other things, that each $n$-extension of $X$ by $Y$ is equivalent to the pullback sequence
$$\xymatrix{0 \ar[r] & Y \ar[r] & I_1 \ar[r] & \cdots \ar[r] & I_n \ar[r] & Q_n \ar[r] & 0 \\
	 		0 \ar[r] & Y \ar[r] \ar@{=}[u] & I_1 \ar[r] \ar[u] &  \cdots \ar[r] & \PB \ar[r] \ar[u] & X \ar[r] \ar[u]_-u &}
$$
for some $u\in \OP L(X,Q_n)$, where $Q_n=\Im(I_n\to I_{n+1})$. Moreover, $u$ has a lifting to $I_n$ if and only if the lower extension is null in $\Ext^n(X,Y)$.\medskip
	
Thus, given an $n$-extension of $X$ by $Y$, and operators $u\in\frak{L}(X,Q_n)$ and $v\in\frak L(K_n,Y)$, their corresponding classes in $\Ext^n(X,Y), \mathcal R^n\frak L(-,Y)X$ and $\mathcal R^n\frak L(X,-)Y$, ``agree'' via the isomorphisms of
(\ref{eq:RnL=Extn}) if and only if there is a commutative diagram	
$$
\xymatrix{
0\ar[r] &K_n \ar[r]\ar[d]_{v} & P_n \ar[r] \ar[d] &\cdots \ar[r] & P_1 \ar[r] \ar[d] & X\ar[r] \ar@{=}[d] & 0\\
0\ar[r] &Y \ar[r]\ar@{=}[d] & Z_n \ar[r] \ar[d] &\cdots \ar[r] & Z_1 \ar[r] \ar[d] & X\ar[r] \ar[d]^u & 0\\
0\ar[r] & Y \ar[r] & I_1 \ar[r]  &\cdots \ar[r] & I_n \ar[r] & Q_n \ar[r] & 0
}
$$	
Also, now without any reference to the middle $n$-extension, the natural isomorphism $\mathcal R^n\frak L(-,Y)X= \mathcal R^n\frak L(X,-)Y$ matches the (respective) classes of $u$ and $v$ precisely when they appear in a commutative diagram
$$
\xymatrix{
0\ar[r] &K_n \ar[r]\ar[d]_{v} & P_n \ar[r] \ar[d] &\cdots \ar[r] & P_1 \ar[r] \ar[d] & X\ar[r] \ar[d]^u & 0\\
0\ar[r] & Y \ar[r] & I_1 \ar[r]  &\cdots \ar[r] & I_n \ar[r] & Q_n \ar[r] & 0
}
$$
As a side effect we have that  $v$ is extensible to $P_n$ if and only if $u$ is liftable to $I_n$.

\section{The derived functors of $\mathfrak K$}\label{sec:RnK}
In this section we focus on the objects associated to the right derived functors of $\mathfrak K(-,Y)$ (via projectives) and $\mathfrak K(X,-)$ (via injectives).

\subsection{$\mathcal R^n\mathfrak K(-,Y)$ via projectives}
Fix a Banach space $Y$ and consider the functor $\fK(-,Y)$.
Given a Banach space $X$, to obtain $\mathcal R^n\mathfrak K(-,Y)X$, one attaches a projective resolution like (\ref{eq:pro-res-X}) to $X$, applies $\mathfrak K(-,Y)$ to $
\xymatrixcolsep{2pc}\xymatrix{
\dots \ar[r]^{d_2} & P_2 \ar[r]^{d_1} & P_1 \ar[r] & 0
}$ to get the complex
$$
\xymatrixcolsep{2pc}
\xymatrix{
0 \ar[r] & \mathfrak K(P_1,Y) \ar[r] &\dots \ar[r] & \mathfrak K(P_{n},Y) \ar[r]^-{d_{n}^*}  & \boxed{\mathfrak K(P_{n+1},Y)}  \ar[r]^-{d_{n+1}^*}  & \mathfrak K(P_{n+2},Y) \ar[r] &\dots
}
$$
and computes the (co)homology at the boxed spot:
$$
\mathcal R^n\mathfrak K(-,Y)X =\dfrac{\ker d_{n+1}^* }{\Im d_{n}^*} = \dfrac{\fK(K_n, Y)}{\fK\text{-extensible}},
$$
with $K_n=\ker d_n=\Im d_{n+1}$.
Therefore, $
\mathcal R^n\mathfrak K(-,Y)X$ can be seen as the quotient of the space of compact operators $K_n\to Y$ by the linear subspace of those operators that admit a compact extension $P_n\to Y$.

\subsection{$\mathcal R^n\mathfrak K(X,-)$ via injectives}
The computation of the derived functors of $\mathfrak K(X, -)$ proceeds by categorical duality: one fixes an injective resolution of $Y$ like (\ref{eq:inj-res-X}), applies $\mathfrak K(X, -)$ to $
\xymatrixcolsep{2pc}\xymatrix{0\ar[r] & I_1\ar[r]^{d^1} & I_2\ar[r]^{d^2} &\dots
}$
 and takes the cohomology of the  resulting complex
$$
\xymatrix{
0\ar[r] & \mathfrak K(X,I_1)\ar[r] & \dots \ar[r] & \mathfrak K(X, I_{n}) \ar[r]^-{d^{n}_*}  & \boxed{\mathfrak K(X, I_{n+1})}  \ar[r]^-{d^{n+1}_*}  & \mathfrak K(X, I_{n+2}) \ar[r] &\dots
}
$$
so that
$$\RK{n}(X,-)Y =  \dfrac{\ker d^{n+1}_* }{\Im d^{n}_*} = \dfrac{\fK(X, Q_n)}{\fK\text{-liftable}},$$
with $Q_n=\Im d^n = \ker d^{n+1}$;
that is, $\RK{n}(X,-)Y$ can be regarded as the quotient of the space of compact operators $X\to Q_n$ by the subspace of those operators that admit a compact lifting to $I_n$.

\subsection{Equivalence between injective and projective derivations}

Our purpose is to show:
\begin{prop}\label{prop:natural} For all Banach spaces $X, Y$ one has
\begin{equation}\label{eq:Ext=R}
\RK{n}(-,Y)X = \RK{n}(X,-)Y
\end{equation}
up to a natural isomorphism that identifies (the equivalence classes of)  $\sigma$ and $\tau$ if and only if there is a commutative diagram
\begin{equation}\label{diag:isom}
		\begin{tikzcd}
			0 \arrow[r] & K_n \arrow[r] \arrow[d, "\tau"] & P_n \arrow[r] \arrow[d] & \cdots \arrow[r] & P_1 \arrow[r] \arrow[d] & X \arrow[r] \arrow[d, "\sigma"] & 0 \\
			0 \arrow[r] & Y \arrow[r]                       & I_1 \arrow[r]                                                    & \cdots \arrow[r]           & I_n \arrow[r]                        & Q_n \arrow[r]                       & 0
		\end{tikzcd}
	\end{equation}
in which all the descending arrows are compact operators.
\end{prop}

\begin{proof} We prove the existence of the required isomorphisms and leave to the reader the verification that they are natural. We proceed by induction on $n\geq1$ and, where possible, we follow the standard proofs that $\RH{n}(-,Y)X = \RH{n}(X, - )Y$ for modules, see \cite[Chapter IV, \S~8]{HS}.
The case $n=1$ appears in \cite[Theorem 3]{mediterranean}.

Let
$
\xymatrixcolsep{1.5pc} \xymatrix{
0\ar[r] & K\ar[r]^{\varkappa}  &P\ar[r]^{\varpi} & X\ar[r] & 0}
$
and
$\xymatrixcolsep{1.5pc} \xymatrix{
0\ar[r] & Y\ar[r]^{\jmath}  &I\ar[r]^{\pi} & Q\ar[r] & 0
}
$
be a projective presentation of $X$ and an injective presentation of $Y$, respectively. We shall construct a linear isomorphism between the spaces
$$
\frac{\fK(K,Y)}{\fK\text{-extensible}}\qquad\text{and}\qquad \frac{\fK(X,Q)}{\fK\text{-liftable}}.
$$
Pick $\tau\in\fK(K,Y)$ and let $\tau_1$ be a compact extension of $\jmath \tau$ obtained  by the first part of Proposition~\ref{prop:K-is-exact}. Since $\pi\tau_1$ vanishes on $K$ one has a commutative diagram
\begin{equation}\label{eq:diag-they-fit}
\xymatrix{
0\ar[r] & K\ar[r]^{\varkappa} \ar[d]^\tau &P \ar[r]^{\varpi} \ar[d]^{\tau_1} & X\ar[r]\ar[d]^\sigma & 0 \\
0\ar[r] & Y\ar[r]^{\jmath}  &I\ar[r]^\pi & Q\ar[r] & 0
}
\end{equation}
with $\sigma\in\fK$. A few points require some checking:
\begin{itemize}
\item The class of $\sigma$ depends only on $\tau$ because a different compact extension $\tau_1'$ would lead to a $\sigma'\in \fK(X, Q)$ with $\sigma-\sigma'$ liftable.
\item For the same reason, the class of $\sigma$ depends linearly on $\tau$, that is, we have a linear map $\fK(K,Y)\to \RK{1}(X, - )Y$...
\item ... whose kernel consists of those operators that admit a compact extension to $P$: if $\lambda\in \fK(X,I)$ is a lifting of $\sigma$, then $\lambda\,\varpi$ is an extension of $\tau$; if $\tau$ turns out to be $\fK$-extensible, one may set $\sigma=0$.
\item In this way the linear map that appeared two-items-ago yields a linear injection $\RK{1}(-,Y)X \to \RK{1}(X, - )Y$.
\item To check surjectivity, dualize: start with $\sigma$ and construct $\tau$, using the second part of Proposition~\ref{prop:K-is-exact}.
\end{itemize}

To complete the proof, let us see that if one couples the classes of $\tau$ and $\sigma$ if there is a commutative diagram of the form \eqref{diag:isom}
 in which all the descending arrows are compact operators, one  gets a linear isomorphism between $\RK{n}(-,Y)X$ and $ \RK{n}(X, - )Y$. 
  We just need a bit extra of detail to also include the kernels:

\begin{equation}
\xymatrix{
0 \ar[r] & K_n \ar[r] \ar[dd]^{\tau} & P_n \ar[rr] \ar[dd]^{\tau_1} \ar[rd] & & P_{n-1} \ar[r] \ar[dd]^{\tau_2} & \cdots \ar[r] & P_1 \ar[r] \ar[dd]^{\tau_n} & X \ar[r] \ar[dd]^{\sigma} & 0 \\
& & & K_{n-1} \ar[dd]^(.3){\tau_{3/2}} \ar[ur] \\
			0 \ar[r] & Y \ar[r]     & I_1 \ar[rr] \ar[rd] &  & I_2 \ar[r]                                                    & \cdots \ar[r]           & I_n \ar[r]                        & Q_n \ar[r]                       & 0\\
			& & & Q_{1} \ar[ur]
}
\end{equation}
Proceeding inductively, the case $n=1$ has already been done. So here the first chunk sends $\tau$ to $\tau_{3/2}$, as in the case $n=1$, while the second sends $\tau_{3/2}$ to $\sigma$ and everything works fine by the induction hypothesis. The splicing is compatible with the equivalences because a compact operator sitting on the $\tau_{3/2}$ spot is $\fK$-liftable (to $I_1$) if and only if it is $\fK$-extensible (to $P_{n-1}$).\end{proof}

As a byproduct of the preceding proof we see that, with the same notations, for each $0<k<n$ one also has natural isomorphisms
\begin{equation}\label{eq:reductRn}
\RK{k}(-,Y)K_{n-k}= \RK{n}(-,Y)X = \RK{n}(X, -)Y =
 \RK{n-k}(X, -)Q_k.
\end{equation}


\subsection{Functoriality}\label{sec:func}
Let's take a look at the action of our functors on arrows (i.e. operators), focusing on the contravariant case $\RK{n}(-,Y)$. Let $X, X'$ be Banach spaces with projective resolutions labeled $(P_k)$ and $(P_k')$, respectively. Given an operator $v:X'\to X$, the lifting property of $P_k'$ provides a commutative diagram
\begin{equation}
		\notag
		\begin{tikzcd}
			0 \arrow[r] & K_n' \arrow[r] \arrow[d, "u"] & P_n' \arrow[r] \arrow[d] & \cdots \arrow[r] & P_1' \arrow[r] \arrow[d] & X' \arrow[r] \arrow[d, "v"] & 0 \\
			0 \arrow[r] & K_n \arrow[r]                       & P_n \arrow[r]                                                    & \cdots \arrow[r]           & P_1 \arrow[r]                        & X \arrow[r]                       & 0
		\end{tikzcd}
	\end{equation}
in which $u$ is well-defined modulo extensible operators and
$$
\RK{n}(-,Y)u : \RK{n}(-,Y)X\to \RK{n}(-,Y)X'
$$
just sends (the class of) $\tau$ to (the class of) $\tau v$. On the other hand, it is simple that if we fix $X$, the assignment $Y \rightsquigarrow  \RK{n}(-,Y)X$ defines a covariant functor whose action on arrows $Y\to Y'$ is given by plain composition. Applying the bifunctor lemma one obtains that $(X,Y) \rightsquigarrow  \RK{n}(-,Y)X$ defines a bifunctor, contravariant on $X$ and covariant on $Y$. In a similar (dual) way  $(X,Y) \rightsquigarrow  \RK{n}(X,-)Y$ is also a bifunctor. The adjective ``natural'' in the statement of Proposition~\ref{prop:natural} refers to these functors.

From now on we write
\begin{equation}\label{eq:Ext=R}
\Ext^n_{\mathfrak K}(X,Y) = \RK{n}(-,Y)X = \RK{n}(X,-)Y,
\end{equation}
and we use the representation that is most convenient at the time. For instance, the chain of isomorphisms in \eqref{eq:reductRn} now becomes the ``reduction formul\ae''
$$
 \Ext^k_{\mathfrak K}(K_{n-k},Y) = \Ext^n_{\mathfrak K}(X,Y) = \Ext^{n-k}_{\mathfrak K}(X,Q_k).
$$

The ideal property of $\fK$ in $\mathfrak{L}$ implies the existence of certain ``Yoneda products'': if $X, Y, Z$ are Banach spaces and $n,m\geq 1$, one has
natural bilinear maps
\begin{itemize}
\item[(a)] $
\Ext^n(X,Z)\times \Ext^m_\fK(Z, Y)\To \Ext_\fK^{n+m}(X,Y)
$;
\item[(b)]$
\Ext^n(X,Y)\times \Ext_\fK^m(Y,Z)\To \Ext_\fK^{n+m}(X,Z)
$.
\end{itemize}

Let us briefly describe the first of these. Consider a projective resolution of $X$ as in \eqref{eq:pro-res-X} and set $K_n=\ker(P_n\to P_{n-1})$. We also need a projective resolution of $Z$ whose elements are decorated with primes.

Pick $u\in\mathfrak{L}(K_n, Z)$ and $\tau\in\fK(K_m', Y)$ and look at
$$
\xymatrix{
K_{n+m}\ar[r] \ar@{..>}[d]_v & P_{n+m}\ar[r]  \ar@{..>}[d]^{u_m}& \dots \ar[r] & P_{n+1}\ar[r] \ar@{..>}[d]^{u_1} & K_n \ar[d]^u\\
K_{m}'\ar[r]\ar[d]_\tau & P_{m}'\ar[r] & \dots \ar[r] & P_{1}'\ar[r] & Z \\
Y}
$$
The dotted arrows
are constructed from right to left, using the lifting property
of the $P_k$ and $v$ is the ``restriction'' of $u_m$ to $K_{n+m}$. The map in (a) takes the class of $(u,\tau)$ in $\Ext^{n}(X,Z)\times \Ext_\fK^{m}(Z,Y)$ to the class of $\tau\,v$ in $\Ext_\fK^{n+m}(X,Y)$.

\subsection{Homology sequences}

The standard brass tacks of homology (cf. Frerick--Sieg \cite[Theorems 5.59 and 5.60]{frer-sieg}) yield the existence of long homology sequences for the functors  $\Ext_{\OP K}^n(E,-)$ and $\Ext_{\OP K}^n(-,E)$. Precisely:

\begin{prop} \label{comp_long} Given a short exact sequence $\xymatrix @C=1.6pc {0\ar[r] &Y \ar[r]^{\imath} & Z \ar[r]^{\pi} & X \ar[r]& 0}$ and a Banach space $E$ there exist exact sequences
		\begin{equation}\label{eq:contra-seq}
			 \begin{tikzcd}[column sep=1cm, row sep=1cm]
			0 \arrow[r] & \mathfrak K(X,E) \arrow[r, "\pi^*"] & \mathfrak K(Z,E) \arrow[r, "\iota^*"]
			\arrow[dl, phantom, ""{coordinate, name=Z}]
			& \mathfrak K(Y,E) \arrow[dll,
			"\omega^0" description, rounded corners=6pt,
			to path={ -- ([xshift=6.5ex]\tikztostart.east)
				|- (Z) [very near end] \tikztonodes
				-| ([xshift=-3ex]\tikztotarget.west)
				-- (\tikztotarget)}] \\
			& \Ext_\mathfrak{K}(X,E) \arrow[r, "\pi^*"] & \Ext_\mathfrak{K}(Z,E) \arrow[r, "\iota^*"]
			\arrow[dl, phantom, ""{coordinate, name=Z}]
			& \Ext_\mathfrak{K}(Y,E) \arrow[dll,
			"\omega^1" description, rounded corners=6pt,
			to path={ -- ([xshift=3ex]\tikztostart.east)
				|- (Z) [very near end] \tikztonodes
				-| ([xshift=-3ex]\tikztotarget.west)
				-- (\tikztotarget)}] \\
			& \Ext_\mathfrak{K}^2(X,E) \arrow[r, "\pi^*"] & \Ext^2_\mathfrak{K}(Z,E) \arrow[r, "\iota^*"]         \arrow[dl, phantom, ""{coordinate, name=Z}] &\Ext^2_\mathfrak{K}(Y,E) \arrow[dll,
			"\omega^2" description, rounded corners=6pt,
			to path={ -- ([xshift=3ex]\tikztostart.east)
				|- (Z) [very near end] \tikztonodes
				-| ([xshift=-2ex]\tikztotarget.west)
				-- (\tikztotarget)}] \\
			& \Ext_\mathfrak{K}^3(X,E)\arrow[r] & \cdots \end{tikzcd}
			\end{equation}	
and
	
		\[ \begin{tikzcd}[column sep=1cm, row sep=1cm]
			0 \arrow[r] & \mathfrak K(E,Y) \arrow[r, "\iota_*"] & \mathfrak K(E,Z) \arrow[r, "\pi_*"]
			\arrow[dl, phantom, ""{coordinate, name=Z}]
			& \mathfrak K(E,X) \arrow[dll,
			"\omega_0" description, rounded corners=6pt,
			to path={ -- ([xshift=6.5ex]\tikztostart.east)
				|- (Z) [very near end] \tikztonodes
				-| ([xshift=-3ex]\tikztotarget.west)
				-- (\tikztotarget)}] \\
			& \Ext_\mathfrak{K}(E,Y) \arrow[r, "\iota_*"] & \Ext_\mathfrak{K}(E,Z) \arrow[r, "\pi_*"]
			\arrow[dl, phantom, ""{coordinate, name=Z}]
			& \Ext_\mathfrak{K}(E,X) \arrow[dll,
			"\omega_1" description, rounded corners=6pt,
			to path={ -- ([xshift=3ex]\tikztostart.east)
				|- (Z) [very near end] \tikztonodes
				-| ([xshift=-3ex]\tikztotarget.west)
				-- (\tikztotarget)}] \\
			& \Ext_\mathfrak{K}^2(E,Y) \arrow[r, "\iota_*"] & \Ext^2_\mathfrak{K}(E,Z) \arrow[r, "\pi_*"]         \arrow[dl, phantom, ""{coordinate, name=Z}] &\Ext^2_\mathfrak{K}(E, X) \arrow[dll,
			"\omega_2" description, rounded corners=6pt,
			to path={ -- ([xshift=3ex]\tikztostart.east)
				|- (Z) [very near end] \tikztonodes
				-| ([xshift=-2ex]\tikztotarget.west)
				-- (\tikztotarget)}] \\
			& \Ext_\mathfrak{K}^3(E,Y)\arrow[r] & \cdots \end{tikzcd}
			\]
			Both exact sequences depend naturally on the starting extension.
	\end{prop}

Note that in the first sequence one identifies $\Ext^n_\mathfrak{K}(-,E)$ with $\RK{n}(-,E)$, while in the second
one identifies $\Ext^n_\mathfrak{K}(E,-)$ with $\RK{n}(E,-)$. Although we are convinced that most readers and authors could write the proof by themselves, we include a reasonably complete construction of the contravariant sequence because actually ``our'' sequence is not exactly the same one obtains (for $\OP L$) in \cite[Chapter IV]{HS}.

\begin{proof} Consider the projective resolution \begin{equation}\tag{$\mathscr P$}
\xymatrixrowsep{1.5pc}
\xymatrix{
\dots \ar[rr] & &  P_3 \ar[rr]^{d_2} \ar[rd]^{q_2}  & &  P_2 \ar[rr]^{d_1}\ar[rd]^{q_1} & &  P_1 \ar[r]^{q_0} & X \ar[r] & 0 \\
&& &K_2 \ar[ur]^{j_1} &&K_1 \ar[ur]^{j_0}&&&
}
\end{equation} of $X$ used in the construction of $\Ext^n_{\OP K}(X,E)=\RK{n}(-,E)X$
and let us also consider a projective resolution of $Y$, written with primes. Following \cite[Proposition~5.43]{frer-sieg} we insert a projective resolution of $Z$ in the middle of the following ``horseshoe'' diagram
	\begin{equation}\label{square}
	\xymatrixcolsep{3pc}
		\xymatrix {            & \vdots \ar[d]^{d_2'} & \vdots \ar[d]^{d_2''}  & \vdots \ar[d]^{d_2}          &   \\
			0 \ar[r] & P_2' \ar[r]^-{\imath_2} \ar[d]^{d_1'} & P_2' \oplus P_2 \ar[r]^-{\pi_2} \ar[d]^{d_1''} & P_2 \ar[r] \ar[d]^{d_1} & 0 \\
			0 \ar[r] & P_1' \ar[r]^-{\imath_1} \ar[d]^{q_0'} & P_1' \oplus P_1 \ar[r]^-{\pi_1} \ar[d]^{q_0''} & P_1 \ar[r] \ar[d]^{q_0} & 0 \\
			0 \ar[r] & Y \ar[r]^-{\imath} \ar[d]                  & Z \ar[r]^-\pi \ar[d]                           & X \ar[r] \ar[d]              & 0 \\
			& 0                                           & 0                                                    & 0                              &
		}
	\end{equation}
making it commutative. 	Here, $\imath_n$ is the inclusion on the first coordinate of $P_n'\oplus P_n$ and $\pi_n$ is the projection on the second; in particular the rows are all exact. The diagram is completed from bottom to top as follows. Let $\widehat q_0: P_1\to Z$ be a lifting of $q_0$ and define $q_0''(a',a)= \imath q_0'(a')+ \widehat{q}_0(a)$. Clearly, $q_0''\imath_1=\imath q_0'$ and $\pi q_0''=q_0\pi_1$, and thus the two lower levels form a commutative diagram.
But $q_0''$ is surjective and, letting $K_1''=\ker q_0''$, we see that the sequence of kernels
$\xymatrixcolsep{0,5cm}\xymatrix{0 \ar[r] & K_1' \ar[r]                   & K_1'' \ar[r]    & K_1 \ar[r]            & 0}$ (the arrows are the restrictions of $\imath_1$ and $\pi_1$) is again exact	and we may repeat the procedure with
$$
	\xymatrixcolsep{3pc}
\xymatrix{
 P_2' \ar[r] \ar[d]_{q_1'}                  & P_2'\oplus P_2 \ar@{..>}[d]_{q_1''} \ar[r]    & P_2    \ar[d]^{q_1}  \ar@{..>}[ld]^{\widehat q_1}        \\
 K_1' \ar[r]                   & K_1'' \ar[r]    & K_1
}
$$
where $\widehat q_1$ is a lifting of $q_1$ and $q_1''(a',a)= \imath_1 q_1'(a')+ \widehat q_1(a)$. Then, take $d_1''$ to be $q_1''$ followed by the inclusion of $K_1''$ into $P_1'\oplus P_1$  and so on. Since the restriction of each 	$\widehat q_n:P_{n+1}\to K_n''$ to $K_{n+1}$ takes values in $K'_{n}$, with $K_0'$ meaning $Y$, we obtain a sequence of operators $g_n: K_{n+1}\to K_n'$ for $n\geq 0$ and we may form the sequence $$
\xymatrixcolsep{1.45pc}
\xymatrix{
\dots \ar[r] & K_3 \ar[r]^{g_2} & K_2' \ar[r] & K_2''\ar[r] & K_2
\ar[r]^{g_1} & K_1' \ar[r] & K_1''\ar[r] & K_1 \ar[r]^{g_0} & Y \ar[r] & Z\ar[r] & X \ar[r] & 0
}
$$
This sequence is not exact, and neither is its image
under $\fK(-,E)$
$$
\xymatrixcolsep{1.2pc}
\xymatrix{
0 \ar[r] & \fK(X,E) \ar[r] &  \fK(Z,E) \ar[r] &  \fK(Y,E) \ar[r]
& \fK(K_1,E) \ar[r] &  \fK(K_1'',E) \ar[r] &  \fK(K_1',E) \ar[r] &   \dots
}
$$
However, the latter becomes exact after taking quotients by the corresponding subspaces of $\fK$-extensible operators, and so it provides the desired homology sequence in the form
\begin{equation}\label{eq:contra-seq}
			 \begin{tikzcd}[column sep=1cm, row sep=1cm]
			0 \arrow[r] & \mathfrak K(X,E) \arrow[r, "\pi^*"] & \mathfrak K(Z,E) \arrow[r, "\iota^*"]
			\arrow[dl, phantom, ""{coordinate, name=Z}]
			& \mathfrak K(Y,E) \arrow[dll,
			"g_0^*" description, rounded corners=6pt,
			to path={ -- ([xshift=6.5ex]\tikztostart.east)
				|- (Z) [very near end] \tikztonodes
				-| ([xshift=-3ex]\tikztotarget.west)
				-- (\tikztotarget)}] \\
			& \dfrac{\fK(K_1,E)}{\fK\text{-extensible}} \arrow[r, "\pi^*"] & \dfrac{\fK(K_1'',E)}{\fK\text{-extensible}} \arrow[r, "\iota^*"]
			\arrow[dl, phantom, ""{coordinate, name=Z}]
			& \dfrac{\fK(K_1',E)}{\fK\text{-extensible}} \arrow[dll,
			"g_1^*" description, rounded corners=6pt,
			to path={ -- ([xshift=3ex]\tikztostart.east)
				|- (Z) [very near end] \tikztonodes
				-| ([xshift=-3ex]\tikztotarget.west)
				-- (\tikztotarget)}] \\
			& \dfrac{\fK(K_2,E)}{\fK\text{-extensible}} \arrow[r, "\pi^*"] & \dfrac{\fK(K_2'',E)}{\fK\text{-extensible}} \arrow[r, "\iota^*"]         \arrow[dl, phantom, ""{coordinate, name=Z}] &\dfrac{\fK(K_2',E)}{\fK\text{-extensible}}\arrow[dll,
			"g_2^*" description, rounded corners=6pt,
			to path={ -- ([xshift=3ex]\tikztostart.east)
				|- (Z) [very near end] \tikztonodes
				-| ([xshift=-2ex]\tikztotarget.west)
				-- (\tikztotarget)}] \\
			& \dfrac{\fK(K_3,E)}{\fK\text{-extensible}}\arrow[r] & \cdots \end{tikzcd}
			\end{equation}	
Recall that
\[ \mathcal R^n\fK(-, E)X = 
\frac{\ker d_{n+1}^*}{\Im d_{n}^*}= \frac{\mathfrak K (K_n,E)}{j_n^*(\mathfrak K (P_n,E))}
=\frac{\mathfrak K(K_n,E)}{\fK\text{-extensible}},
\]
where $d_0^*=0$, and analogously for $Y$ and $Z$. Also, the arrows between objects in the same level are as expected since
$$
\pi^*_n : \dfrac{\fK(K_n,E)}{\fK\text{-extensible}}\To \dfrac{\fK(K_n'',E)}{\fK\text{-extensible}}
$$
is the image of $\pi: Z\to X$ under $\RK{n}(-,E)$, and the same with $\imath_n$. One must still check exactness, for which we give the details at the $K_n$-th spots; the other cases are easier. The relevant section of the sequence is therefore
$$
\xymatrixrowsep{1pc}
\xymatrix{
\cdots\ar[r] & \dfrac{\fK(K_n',E)}{\fK\text{-extensible}} \ar[r]^{g_n^*} \ar@{=}[d] &
  \dfrac{\fK(K_{n+1},E)}{\fK\text{-extensible}} \ar[r]^{\pi_{n+1}^*}  \ar@{=}[d] &
   \dfrac{\fK(K_{n+1}'',E)}{\fK\text{-extensible}} \ar[r]  \ar@{=}[d] & \cdots\\
   & \mathcal R^n\fK(-, E)Y  & \mathcal R^{n+1}\fK(-, E)X  &\mathcal R^{n+1}\fK(-, E)Z
}
$$
Think of the following part of the $n$-th horsehoe diagram (and ignore the dotted arrows)

$$
\xymatrixcolsep{3pc}\xymatrixrowsep{3pc}\xymatrix
{ K_{n+1}' \ar[r]^{\imath_{n+1}} \ar[d]^{\jmath_n'} &K_{n+1}'' \ar[r]^{\pi_{n+1}} \ar[d]^{\jmath_n''} & K_{n+1} \ar@{..>}[dddrr]^\alpha \ar@{-->}@/_11pc/[ddll]_{g_n} \ar[d]^{\jmath_n} \\ P_{n+1}'\ar[r]^{\imath_{n+1}} \ar[d]^{q_n'} & P_{n+1}'' \ar@{..>}@/_3pc/[ddrrr]^{\beta} \ar[r]^{\pi_{n+1}} \ar[d]^{q_n''} & P_{n+1} \ar[d]^{q_n}
 \\
K_{n}' \ar[r]^{\imath_n}   \ar@{..>}@/_1pc/[drrrr]_{\gamma}             & K_{n}'' \ar[r]^{\pi_n}       & K_{n} \\
&&&& E
}
$$
The containment $\Im g_n^*\subset \ker \pi_{n+1}$  is clear because the composition $g_n\pi_{n+1}: K_{n+1}''\to K_{n+1}\to K_{n}'$ actually extends to $P_{n+1}''=
P_{n+1}'\oplus P_{n+1}$. Indeed, one has $q_n''(u,v)= \imath_nq_n'u + \widehat{q}_n v$; if $(u,v)\in K_{n+1}''$, then   $\imath_nq_n'u + \widehat{q}_n v=0$ and letting $g(u,v)=-q_n'(u)$ we
get an extension of $g_n\pi_{n+1}$. To obtain the reverse containment, pick $\alpha\in \fK(K_{n+1}, E)$ such that $\alpha\pi_{n+1} = 0$ in $\mathcal R^{n+1}\fK(-, E)Z$
and let $\beta\in \fK(P'', E)$ be a witness, so that $\beta \jmath_n''=  \alpha\pi_{n+1}$. By commutativity,
$\beta\imath_{n+1} \jmath_n'=  \alpha\pi_{n+1}\imath_{n+1}=0$.
Hence there exists $\gamma\in\fK(Y, E)$ such that $\gamma q_n'=\beta \imath_{n+1}$. Does this $\gamma$
satisfy $\gamma g_n = \alpha$ in $\mathcal R^{n+1}\fK(-, E)Y$? Not exactly. Recall that $q_n''(u,v)= q_n'(u) + \widehat{q}_n(v)$, so that $q_n'(u) = - \hat{q}_n(v)$ for $(u, v)\in K_{n+1}''$. Now assume $v\in K_{n+1}$ and choose $u\in P_{n+1}'$ such that $(u,v)\in K_{n+1}''$. Then
$$
\gamma g_n(v)=  \gamma \widehat{q}_n(v)= -\gamma q_n'(u)= -\beta(u,0)
= -\beta(u,v) + \beta(0,v) = -\alpha(v) +\beta(0,v).
$$
Hence the compact operator $w\in P_{n+1}\mapsto \beta(0,w)\in E$ shows that  $\gamma g_n = -\alpha$ in $\mathcal R^{n+1}\fK(-, E)Y$ and, therefore, $g_n^*(-\gamma)=\alpha$ in $\mathcal R^{n+1}\fK(-, E)Y$, which is enough.\medskip

An analogous (dual) process using injective resolutions yields the covariant sequence. The natural character of the two sequences is fairly obvious.\end{proof}

\subsection{Other snakes} Let us remark at this point that ``our'' sequence \eqref{eq:contra-seq} is not exactly the same a professional algebraist would have obtained. Indeed, focus on the first connection morphism of the $\Hom$-$\Ext$ sequence and look at the construction Hilton--Stammbach display in \cite[Theorem III.5.3 and Lemma III.5.4]{HS}. Doing that, we should  have started with a diagram of projective presentations (a lot of indices have been omitted)
\begin{equation}\label{other-approachesnet}
		\xymatrix {
0\ar[r]	 & K' \ar[r] \ar[d] & K'' \ar[r] \ar[d] & K \ar[r] \ar[d] & 0 \\
			0 \ar[r] & P' \ar[r] \ar[d]^{q'} & P'' \ar[r]^{\pi_1} \ar[d]^{q''} & P \ar[r] \ar[d]^{q} & 0 \\
			0 \ar[r] & Y \ar[r]^{\imath}                   & Z \ar[r]^{\pi}                            & X \ar[r]               & 0 \\
		}
	\end{equation}
Then apply $\fK(-, E)$  to the first two rows to get
\begin{equation}\label{other-approaches}
				\xymatrix{ 0 \ar[r] & \fK(P,E) \ar[r]  \ar[d] & \fK(P'',E)  \ar[d] \ar[r]
 & \fK(P',E) \ar[r]  \ar[d] & 0
 \\
 0 \ar[r] & \fK(K,E) \ar[r] & \fK(K'',E) \ar[r]
 & \fK(K',E)
		}
	\end{equation}
(the upper row in this diagram is exact because the middle row in the preceding one splits) and finally apply the snake lemma \cite[Lemma III.5.1]{HS} to obtain an operator
$$
\xymatrixcolsep{2.5pc}
			\xymatrixrowsep{1.5pc}
\xymatrix{
\fK(Y,E) \ar@{=}[d] \ar[rr]^{\tilde\omega^0} &  & \Ext_{\fK}(X, E) \ar@{=}[d]^{\text{def}} \\
\ker\big(\fK(P',E)\to \fK(K',E) \big) & & \coker \big(\fK(P,E)\to \fK(K,E) \big)
}
$$
making the sequence
$$
\xymatrixcolsep{1.4pc}
\xymatrix{0 \ar[r] & \fK(X,E) \ar[r] & \fK(Z,E) \ar[r] & \fK(Y,E) \ar[r]^-{\tilde\omega^0} &
\Ext_\fK(X,E) \ar[r] & \Ext_\fK(Z,E) \ar[r] & \Ext_\fK(Y,E)
}
$$
automatically exact. It is easy to check that this connecting morphism $\tilde\omega^0$ is plain composition with the operator $f: K\To Y$ obtained from Diagram (\ref{other-approachesnet}) in the following snakish way: Let $s\in\frak L(P, P'')$ be a section of $\pi$ and consider the projection of $P''$ onto $P'$
given by  $\pi' = \mathbf 1_{P''} - s \pi_1$.
Now,
pick $k\in K$, choose $k''\in K''$ such that $\pi_1(k'')=k$ and set $f(k) = q'\pi'(k'')$. The choice of $k''$ does not affect
$f(k)$ because the restriction of $q'\pi'$ to $K'$ is zero (if $k'\in K'$, then $q'\pi'(k')= q'(k')=0$).
  Moreover, and this is the funny part, since $q''(k'')=0$ it turns out that $0 = q''(k'') = f(k) + \widehat q(k)$ which means that, while our choice for the map $g: K\to Y$ was $g(k)=\widehat q(k)$, the algebraist's choice yields $f(k)$ (see \cite[p.~104]{HS}) and they are therefore opposites: $f+g=0$, and the same occurs with the associated composition maps. This infection propagates, and in general, $\tilde\omega^n =(-1)^{n+1}\omega^n$ for the contravariant sequences. Thus, in a sense, the sequence \eqref{eq:contra-seq} is more natural than the traditional one \cite[Chapter IV, \S 9]{HS} (see the final comments on pp. 154--155). The covariant construction is not affected and $\tilde\omega_n =\omega^n$. This difference may be puzzling, but an explanation is that while our sequence is obtained by categorical duality, using  injectives, that is, $\Ext^n_\mathfrak{K}(E, X)=\RK{n}(E, -)X$, if one wanted to adapt Hilton--Stammbach's Theorem III.5.2 and Proposition IV.7.4, and their proofs, then  $\Ext^n_\mathfrak{K}(E, X)$ should be interpreted as $\RK{n}(-,X)E$.

\section{Comparing homologies}\label{sec:comparing}


Since it is irrefutable that every compact operator is an operator, it is clear that one has natural linear maps
$\eta: \Ext_{\mathfrak K}^n(X,Y)\to \Ext^n(X,Y)$ depending on $X,Y$ and $n$ --- all of which make a natural transformation of functors, of course. Precisely, if we interpret $\Ext_{\mathfrak K}^n(X,Y)$ as $\RK{n}(-,Y)X$ and $\Ext^n(X,Y)$ according to Yoneda, it should be clear that $\eta$ takes (the class of) $\tau\in\fK(K_n,Y)$ into the (class of the) pushout sequence
\begin{equation}
		\notag
		\begin{tikzcd}
			0 \arrow[r] & K_n \arrow[r] \arrow[d, "\tau"] & P_n \arrow[r] \arrow[d] & \cdots \arrow[r] & P_1 \arrow[r]                                & X \arrow[r]                                & 0  \\
			0 \arrow[r] & Y \arrow[r]             & \PO \arrow[r]                                                   & \cdots \arrow[r] & P_1 \arrow[r] \arrow[u, equal] & X \arrow[r] \arrow[u, equal] & 0
		\end{tikzcd}
	\end{equation}

 In this section we compare $\Ext^n(X,Y)$ and $\Ext_{\mathfrak K}^n(X,Y)$ by means of $\eta$,  mainly for  $n=1$. As we shall see, the image of
$\Ext_{\mathfrak K}^n$ in (Yoneda) $\Ext^n$ leads to a quite interesting class of $n$-extensions that we have called  ``compact $n$-extensions''.

Most of what is known about $\Ext^n(X,Y)$ for classical Banach spaces $X,Y$ and $n=1$ appears in \cite{hmbst}; for $n\geq 2$, see \cite{wod, ccg}.

\subsection{Compact extensions} We prepare the way towards the definition with a lemma. The case $n=1$ appears in \cite[Theorem 3]{mediterranean}.

\begin{lem}\label{lem:compact-ext-equiv} They are equivalent:

\emph{(a)}\quad There exists $\tau\in\fK$ and a commutative diagram
\begin{equation}
		\notag
		\begin{tikzcd}
			0 \arrow[r] & K_n \arrow[r] \arrow[d, "\tau"] & P_n \arrow[r] \arrow[d] & \cdots \arrow[r] & P_1 \arrow[r] \arrow[d] & X \arrow[r] \arrow[d, equal] & 0 \\
			0 \arrow[r] & Y \arrow[r]                       & Z_n \arrow[r]                                                    & \cdots \arrow[r]           & Z_1 \arrow[r]                        & X \arrow[r]                       & 0
		\end{tikzcd}
	\end{equation}

\emph{(b)} \quad There exists $\sigma\in\fK$ and a commutative diagram
 \begin{equation*}
	 	\label{eq:compact-pb} \begin{tikzcd}
	 		0 \arrow[r] & Y \arrow[r] & I_1 \arrow[r] & \cdots \arrow[r] & I_n \arrow[r] & Q_n \arrow[r] & 0 \\
	 		0 \arrow[r] & Y \arrow[r] \arrow[u, equal] & Z_n \arrow[r] \arrow[u] &  \cdots \arrow[r] & Z_1 \arrow[r] \arrow[u ]& X \arrow[r] \arrow[u, "\sigma"] & 0
	 	\end{tikzcd}
	 \end{equation*}
\end{lem}

\begin{proof} We prove the implication (a)$\implies$(b). The converse follows by (categorical)
duality. Consider the diagram (ignore for the moment the dotted arrows)
\begin{equation}\label{homotopic}
\xymatrix{
0 \ar[r] & K_n \ar[r] \ar[dd]_{\tau} & P_n \ar[r] \ar@{-->}@/_2pc/[dddd]_(.4){\tau_1} \ar[dd]^{u_1}  & \cdots  \ar[r] & P_2 \ar@{-->}@/_2pc/[dddd]_(.4){\tau_{n-1}} \ar[r] \ar[dd]^{u_{n-1}}  & P_{1} \ar@{-->}@/_2pc/[dddd]_(.4){\tau_n} \ar[r] \ar[dd]^{u_n}  & X \ar@{-->}@/^2pc/[dddd]^(.4){\sigma} \ar[r]   \ar@{=}[dd] & 0 \\
&   \\
			0 \ar[r] & Y \ar[r] \ar@{=}[dd]    & Z_n \ar[dd]^{v_1} \ar[r]           & \cdots \ar[r]                         & Z_{2} \ar[dd]^{v_{n-1}} \ar[r]  & Z_1 \ar[dd]^{v_n} \ar[r]                                                     & X \ar[r]   \ar[dd]_{v}                    & 0\\
			&&& & & \\
0 \ar[r] & Y \ar[r]    & I_1 \ar[r]           & \cdots \ar[r]                         & I_{n-1} \ar[r]   & I_n \ar[r]                                                     & Q_n \ar[r]                       & 0\\
}
\end{equation}
The top half of the diagram is taken from the hypothesis (a). The bottom half is obtained from the injective presentation of $Y$. This, however,  does not guarantee the compactness of any $v_k$ or of $v$. The dashed arrows are obtained from the first part of Proposition~\ref{prop:K-is-exact}: compact operators with values on $\mathscr L_\infty$ admit compact extensions anywhere. They are therefore all compact, but do not necessarily commute with arrows beginning or ending in the spaces $Z_k$, nor with $v$. Nevertheless, both
$(\tau, \tau_1, \dots, \tau_n,\sigma)$ and $(\tau, v_1 u_1, \dots, v_n u_n,v)$ are morphisms of complexes and, therefore, $\sigma-v$ is liftable. 
 The morphism required in (b) is the following: let $\lambda: X\to I_n$ be a  lifting for $\sigma-v$, not necessarily compact. Since
$$
\xymatrix{
Z_{2} \ar[d] \ar[r] & Z_{1} \ar[d]^{v_n+\lambda\pi} \ar[r]^\pi & X \ar[d]^{\sigma=v+\sigma-v} \\
I_{n-1} \ar[r] & I_{n} \ar[r] & Q_n
}
$$
remains commutative, set $({\mathbf 1}_Y, v_1, \dots, v_n+\lambda\pi, \sigma)$.\end{proof}

The subjacent property of the ideal $\OP K$ described in the preceding lemma was called in \cite[Definition 7]{mediterranean} ``to be \emph{balanced}".


We are now in position to define \emph{compact} $n$-extensions of $X$ by $Y$ as those exact sequences that satisfy the equivalent conditions of the just proved lemma. One has:

\begin{prop}
An $n$-extension of $X$ by $Y$ is compact if and only if its class in $\Ext^n(X,Y)$ belongs to the image of
$\eta: \Ext_{\mathfrak K}^n(X,Y)\to \Ext^n(X,Y)$
\end{prop}

\begin{proof} The ``only if'' part is trivial.
The converse is clear because if $\mathscr Z$ satisfies the first condition of Lemma~\ref{lem:compact-ext-equiv} and $\mathscr Z \gtrdot \mathscr Z'$, then $\mathscr Z'$ fits in a similar diagram; while if $\mathscr Z$ satisfies the second condition and $\mathscr Z \lessdot \mathscr Z'$, then so does $\mathscr Z'$.
\end{proof}

Of course, one does not expect $\eta$ to be surjective in general, if only because there are many operators that are not compact. Actually, in view of Proposition~\ref{prop:K-is-exact}
if $X$ is an $\mathscr L_1$-space, then
$\Ext_{\mathfrak K}^n(X, Y)=0$ for all Banach spaces $Y$ and every $n\geq 1$, and if
if $Y$ is an $\mathscr L_\infty$-space, then
$\Ext_{\mathfrak K}^n(X, Y)=0$ for all Banach spaces $X$ and every $n\geq 1$. However:
\begin{itemize}
\item There is an $\mathscr L_1$-space $Y$ such that $\Ext^n(L_1, Y)\neq 0$ for each $n\geq 1$; see \cite[pp. 541--542]{ccg}. Of course, $L_1$ itself is an $\mathscr L_1$-space.
\item If $Y$ is an $\mathscr L_\infty$-space which is not injective, then its injective
 presentation does not split and in particular $\Ext(I/Y, Y)\neq 0$. There is a lot of work around trivial and nontrivial exact sequences $0\to C(K) \to Z \to X \to  0$ in which $X$ is a classical Banach space; see \cite[Chapter~8]{hmbst}.
\item Since $\ell_\infty/c_0$ fails to be injective no injective 2-resolution of $c_0$ is zero in $\Ext^2$; see \cite[Theorem 1.25]{aviles} and \cite[Proposition 4.2]{ccg}.
\end{itemize}

\subsection{Compact short exact sequences}
The reader arrived to this point may rightly argue that the above examples are rather trivial: after all, the null map is not usually surjective. To exhibit more interesting examples we focus on the case $n=1$ adding some ideas from intermediate Banach space theory:
\begin{itemize}
\item A Banach space is an ultrasummand if it is complemented in its bidual --- or in any other dual space.
\item A short exact sequence of Banach spaces splits locally if the dual sequence splits.
\item A subspace $Y\subset Z$ is said to be locally complemented if the natural sequence $0\to Y\to Z \to  Z/Y\to 0$ splits locally.
\item A Banach $X$ space has the $\lambda$-AP if for every finite dimensional subspace $F\subset X$ and every $\varepsilon>0$ there is a finite rank endomorphism $f\in\frak L(X)$, with $\|f\|<\lambda+\varepsilon$, such that $f(x)=x$ for all $x\in F$. A Banach space has the bounded approximation property (BAP) if it has the $\lambda$-AP for some $\lambda\geq 1$.
\end{itemize}

The reader is advised that these definitions are not the usual ones, but equivalent formulations that suit our present purposes. The notions just introduced are connected to the subject of this paper by the following facts (\cite[Theorems 3.4 and 3.5]{kalton-CEP}):
\begin{itemize}
\item[$\bigstar$] A short exact sequence that splits locally and has an ultrasummand in the spot of the subspace splits.
\item[$\bigstar$] $Y$ is locally complemented in $Z$ if and only if it has the compact extension property (CEP): every compact operator from $Y$ to any Banach space has a compact extension to $Z$.
\end{itemize}

We also need a technicality that can be considered a particular case of \cite[Lemma~4.3.3]{hmbst}:

\begin{lem}\label{lem:433} Assume one has morphism of sequences
$$\xymatrix{0\ar[r] &Y \ar[d]_\beta \ar[r]& Z \ar[r]\ar[d]& X \ar[r]\ar[d]^\alpha& 0\\
			0\ar[r] & Y\ar[r]& Z' \ar[r]& X \ar[r]& 0
		}$$
in which either $\beta\in\fK$ and $\alpha={\bf I}_X$ or
$\alpha\in\fK$ and $\beta={\bf I}_Y$.
If the rows are equivalent in $\Ext(X,Y)$ then both are trivial.
\end{lem}

\begin{prop} Given a commutative diagram with exact rows
\begin{equation*}
\xymatrix{
(\mathscr Z) &0\ar[r] & Y \ar[r] \ar[d]_\tau & Z \ar[r]\ar[d] & X\ar[r]\ar@{=}[d]&0 &\\
(\mathscr L) &0\ar[r] & Y'  \ar[r] & \mathscr L_\infty  \ar[r] & X \ar[r]&0 &\quad}
\end{equation*}
with $\tau\in\fK$, if $\mathscr L$ is not trivial then $\mathscr Z$ cannot be compact. Given a commutative diagram $$\xymatrix{
(\mathscr L) &0\ar[r] & Y \ar[r]\ar@{=}[d]  & \mathscr L_1 \ar[r] \ar[d] & X'\ar[r]\ar[d]^\sigma &0 &\\
(\mathscr Z) & 0\ar[r] &  Y \ar[r]  & Z \ar[r] & X \ar[r]&0 &\quad}
$$
with $\sigma\in\fK$, if $\mathscr L$ is not trivial then $\mathscr Z$ cannot be compact. In particular: if $X$ is not projective then no projective presentation of $X$ is compact; and if $Y$ is not injective, no injective presentation of $Y$ is compact.
\end{prop}

\begin{proof} Let us prove the first part. The second is analogous. Since $\tau$ is compact, Proposition~\ref{prop:K-is-exact} yields a commutative diagram
$$\xymatrix{
0\ar[r] & Y \ar[r] \ar[d]_\tau & Z \ar[r]\ar[d] & X\ar[r]\ar[d]^{\sigma}&0\\
0\ar[r] & Y'  \ar[r] & \mathscr L_\infty  \ar[r] & X \ar[r]&0}$$
with $\sigma$ compact. Hence
$\mathscr L \sim \tau\mathscr Z\sim \mathscr L \sigma $ with $\sigma\in\fK$ and the preceding lemma yields that $\mathscr L$ is trivial.\end{proof}

Which other exact sequences fail to be compact? Recall that an operator $u:A\to B$ is said to be strictly singular if the restriction to the infinite dimensional subspaces of $A$ is
never an embedding. Short exact sequences with strictly singular quotient map are often called singular. These include the popular Kalton--Peck sequences $0\to \ell_p\to Z_p\to \ell_p\to 0$ and many others, see \cite[Chapter~9]{hmbst}.

\begin{prop}\label{prop:nsingular} No singular short exact sequence is compact.
\end{prop}
\begin{proof}
By \cite[Proposition 9.1.3 (i)]{hmbst} it suffices to prove the following: \emph{Let $B, Y$ be Banach spaces and $K\subset B$ a subspace of infinite codimension. Then, for every compact $\tau:K\to Y$ there is $A\subset B$ containing $K$, with $A/K$ infinite dimensional, and a compact $\tau':A\to Y$ extending $\tau$.}

According to Diestel \cite[p. 7]{diestel} it was Grothendieck who first proved that every compact operator factors ``compactly'' through a subspace of $c_0$. More precisely: if $\tau:K\to Y$ is compact, then there exists $H\subset c_0$ and compact operators $\alpha:K\to H$ and $\beta:H\to Y$ such that $\tau=\beta\alpha$. Thus, letting $X=B/K$, we have a commutative diagram
$$\xymatrix{
0\ar[r] & K \ar@/_2pc/[dd]_(.4){\tau} \ar[r] \ar[d]^\alpha & B \ar[d]^{\alpha'}\ar[r]^\pi  & X  \ar[r]&0\\
&H\ar[r] \ar[d]^\beta & c_0\\
& Y }$$
with $\alpha'\in \OP K$. Let $(x_n)_{n\geq 1}$ be a normalized basic sequence in $X$ (such a thing exists insofar as $X$ is of infinite dimension, see \cite[Corollary 1.5.3]{a-k}). Extending the coordinate functionals by Hahn--Banach one obtains a bounded sequence $(x_n^*)_{n\geq 1}$ in $X^*$ such that $x_n^*(x_k)=\delta_{n,k}$ for $n,k\in\mathbb N$. Let $(b_n)_{n\geq 1}$ be a bounded sequence in $B$ such that $\pi b_n = x_n$. Since $\alpha'$ is compact there is $z\in c_0$ and an infinite set $M\subset \mathbb N$ such that $\|\alpha'b_{m} - z\|\leq 2^{-m}$ for all $m\in M$. We form the following nuclear operator $\nu: X\to c_0$
$$\nu = \sum_{m\in M} x_m^* \otimes \big(\alpha'b_{m} - z\big).
$$
Observe that for every $m\in M$, one has
$
\alpha'b_{m} - \nu \pi b_{m}   = \alpha'b_{m} - \nu x_{m}   = z$
and therefore the restriction of $\alpha' - \nu\pi$ to $A$, the subspace spanned by $\{b_m: m\in M\}$ in $B$, takes values in $L$, the subspace spanned by $z$ in $c_0$. Pick now a bounded linear projection $e: H + L \to H$. Then the restriction of $\beta e ( \alpha' - \nu\pi)$ to $K+A$ is an extension of $\tau$ and since $K$ has infinite codimension in $\overline{K+A}$ we are done.\end{proof}

We now
show that $\eta: \Ext_{\mathfrak K}(X,Y)\to \Ext(X,Y)$
is not necessarily injective,
something that can be of interest for Banach spacers in the following formulation:
	
\begin{ex}\label{ex:posno}
A compact operator between Banach spaces that admits an extension to some superspace does not necessarily admit a compact extension.
\end{ex}

	\begin{proof} We start with a sequence $0\to Y\to Z\to X\to 0$ and an operator $\alpha: Y \to E$ such that
the lower sequence in the pushout diagram	
		$$\xymatrix{0\ar[r] &Y \ar[d]_\alpha \ar[r]& Z \ar[r]^\pi\ar[d]_{\alpha'}& X \ar[r]\ar@{=}[d]& 0\\
			0\ar[r] & E\ar[r]^-\imath& \PO \ar[r]^{\pi'}& X \ar[r]& 0
		}$$
does not split.
Obviously, $\alpha':Z\to \PO$ is an extension of $\imath\alpha$.
However, this $\alpha'$ is not compact, regardless of what $\alpha$ is. Assume that $\alpha$ is compact and that some compact extension $\beta: Z\to \PO$ of $\imath \alpha$ exists. That generates a commutative diagram
		$$\xymatrix{0\ar[r] &Y \ar[d]_\alpha \ar[r]& Z \ar[r]^{\pi}\ar[d]_{\beta}& X \ar[r]\ar[d]_{\gamma}& 0\\
			0\ar[r] & E\ar[r]& \PO \ar[r]^{Q}& X \ar[r]& 0}$$
If $\mathscr Z$ denotes the starting sequence, then the pushout sequence is $\alpha\mathscr Z$ and the preceding diagram shows that $\alpha\mathscr Z\sim (\alpha\mathscr Z)\gamma$.
But  $\gamma$ has to be compact since $\gamma \pi = \pi'\beta$ is compact. Applying again Lemma~\ref{lem:433} we get $\alpha \mathscr Z\sim 0$ in flagrant contradiction with ourselves.
	\end{proof}

Let us show how the existence of a compact operator $\alpha:Y\to E$ that is extensible but not $\fK$-extensible to $Z$ leads to a (nonzero) compact operator in the kernel of $\eta: \Ext_{\mathfrak K}(X,E)\to \Ext(X,E)$, where $X=Z/Y$.

Take a projective presentation of $X$ and look at the diagram (we have omitted arrows beginning or ending in zeros)
$$\xymatrix{K \ar[d]_u \ar[r] & P \ar[r] \ar[d] & X \ar@{=}[d]\\
Y \ar[d]_\alpha \ar[r]& Z \ar[r]\ar[d]& X \ar@{=}[d]\\
			 E\ar[r]& \PO \ar[r] & X }
$$
Clearly,  if $\alpha$ extends to $Z$,  the composition $\alpha u$ has an extension to $P$. On the other hand, the upper left square is a pushout, so if $\alpha u$ had  an extension   $\tilde\alpha\in\fK(P, E)$ then since $\alpha u= \tilde\alpha$ on $K$ the induced operator $Z\to E$ (which clearly extends $\alpha$) should be compact since $Y\oplus P\to Z$ is surjective and $\alpha\oplus \tilde\alpha : Y\oplus P\to E$ is compact.

\medskip

There is another way of reading the example: if one has a compact operator that cannot be extended to some superspace then there is a compact operator admitting an extension but not a compact extension. Compact operators that cannot be extended occur quite often: if $Y$ is not an $\mathscr L_\infty$-space but has a basis, then there is an enlargement $Z\supset Y$ and a compact endomorphism of $Y$ that cannot be extended to an operator $Z\to Y^{**}$, see Lindenstrauss \cite[Theorem 3.4]{memoir}. More on this in Section~\ref{sec:misc}.

\medskip

	We consider now whether $\Ext_{\mathfrak K}(X,E)=0$ implies $\Ext(X,E)=0$ or viceversa.
	
	\begin{thm} \label{Ext-K} Let $X$ and $Y$ be Banach spaces.
		\begin{enumerate}
			\item[(a)] If $Y$ has the BAP and $\Ext(X, Y)=0 $ then $\Ext_{\mathfrak K}(X,Y)=0$.
			\item[(b)] If $Y$ is an ultrasummand, either $Y$ or $X$ has the BAP and $\Ext_{\mathfrak K}(X,Y)=0$, then $\Ext (X,Y)=0$.
		\end{enumerate}
	\end{thm}
	\begin{proof} (a) By the open mapping theorem, if  $\Ext(X,E)=0$ then for every projective presentation
$0\to K\to P\to X\to 0$ there is a constant $C$ such that every operator $u:K\to Y$ admits an extension $\tilde{u}:P\to Y$ with $\|\tilde{u}\|\leq C\|u\|$. If $Y$ has the $\lambda$-AP and $u$ has finite rank, we may take a finite rank $v\in \mathfrak L(Y)$, with $\|v\|\leq \lambda+\varepsilon$ that is the identity on the range of $u$. Clearly, $v\tilde{u}$ is a finite rank extension of $u$, with $\|v\tilde{u}\|\leq (\lambda+\varepsilon) C\|u\|$. It follows that the restriction map $\frak F(P,Y)\to \frak F(K,Y)$ is open and since
$\frak F(K,Y)$ is dense in $\frak K(K,Y)$ it follows that 	$\frak K(P,Y)\to \frak K(K,Y)$ is onto --- and open.

\smallskip

(b) Let $0\to K\to P\to X\to 0$ be a projective presentation of $X$. If $\Ext_{\mathfrak K}(X,Y)=0$, then there is a constant $C$ such that every compact operator $u: K\to Y$ admits a (compact) extension $\tilde{u}:P\to Y$ with $\|  \tilde{u}\|\leq C\|u\|$. If $Y$ has the $\lambda$-AP, for every finite dimensional $F\subset Y$ we may choose a finite rank (in particular compact) $f_F\in \frak L(Y)$ with $\|f_F\|\leq\lambda+\varepsilon$ and $f_F(y)=y$ for all $y\in F$. Let $\delta:Y\to Y^{**}$ be the natural inclusion map and let $\pi: Y^{**}\to Y$ be a bounded projection. Now, pick $u\in \frak L(K,Y)$ and for each finite dimensional $F\subset Y$ let $u_F:K\to Y$ be an extension of $f_F u$ with $\|u_F\|\leq C\|f_F u\|\leq C(\lambda+\varepsilon)\|u\|$.

Let $\mathscr U$ be an ultrafilter on the set of all finite dimensional subspaces of $Y$ that refines the order filter and define $\tilde{u}:P\to Y^{**}$ by
$$
\tilde{u}(x) = \lim_{\mathscr U(F)} \delta u_F(x),
$$
where the limit is taken in the weak* topopology of $Y^{**}$. Then $
\|\tilde{u}\|\leq C(\lambda+\varepsilon)\|u\|$) and, clearly, $\tilde{u}|_K= \delta u$, so that $\pi \tilde{u}$ is the required extension of $u$.

\smallskip

If, instead, we assume that $X$ has the BAP then so does $K$, by a result of Lusky (see \cite[Proposition 5.3.3]{hmbst}) and we may take, for each finite dimensional $F\subset K$ a finite rank operator $g_F\in \frak L(K)$ such that $g_F(x)=x$ for $x\in F$, with $\sup_F\|g_F\|<\infty$. Then we can consider $ug_F$ and   proceed as before.
	\end{proof}

From the preceding theorem and the corresponding counterparts for $\Ext$ (see \cite[Proposition 5.2.20]{hmbst}) we get:
	
\begin{cor}
$\Ext_{\fK}(\ell_p, \ell_q)\neq 0$ for every $p\in(1,\infty], q\in[1,\infty)$, and the same if $\ell_\infty$ is replaced by $c_0$.
\end{cor}

The hypothesis ``$Y$ is an ultrasummand" in the second part of the theorem above is necessary:
if $Y$ is an $\mathscr L_\infty$-space, in particular $C[0,1]$ or $c_0$, then $\Ext_{\mathfrak K}(-, Y)=0$ by Proposition~\ref{prop:K-is-exact}; however $\Ext(X, C[0,1])\neq 0$ if $X$ is $c_0$ or $\ell_p$ for $p\in(1,\infty]$  and $\Ext(X,c_0)\neq 0$ if $X$ is a nonseparable WCG Banach space, see \cite[Proposition 8.6.4 and 4.5.13]{hmbst}.

These examples also show that assertion (a) in Theorem \ref{Ext-K} cannot be reversed. We do not know whether some hypothesis on the approximation property is really necessary there, but note that one can replace the BAP by the (weaker) bounded compact approximation property, see \cite[Section~8]{casazza}.

\subsection{Nonlinear characterization of compact extensions}\label{nonlinear}

So far, we have done our best to avoid using quasilinear maps, but we can no longer do without them. We refer the reader to the third chapter of \cite{hmbst} for an account of quasilinear deeds in exact sequence affairs. A map $\Phi: X\to Y$, acting between normed spaces, is said to be $1$-linear if it
is homogeneous and there is a constant $Q_1$ such that for every $n\geq 2$ and every $x_1,\dots,x_n\in X$ one has
$$
\left\| \Phi\left(\sum_{i\leq n} x_i\right)- \sum_{i\leq n} \Phi(x_i )\right\| \leq Q_1 \sum_i\|x_i\|.
$$
The least constant for which the inequality holds (for all $n$) shall be referred to as the 1-linearity constant of $\Phi$, denoted by $Q_1(\Phi)$.
If the preceding inequality holds merely for $n=2$ we say that $\Phi$ is quasilinear.\medskip

Given a quasilinear map $\Phi$ one can endow the product space $Y\times X$ with the quasinorm
$\|(y,x)\|_\Phi= \|y-\Phi(x)\|+\|x\|$, which is equivalent to a norm if and only if $\Phi$ is 1-linear. If $Y\oplus_\Phi X$ denotes the resulting quasinormed space
one has an exact sequence
$$
\xymatrix{
		0 \ar[r] & Y \ar[r]^-\imath & Y\oplus_\Phi X \ar[r]^-\pi & X \ar[r] & 0}
$$
where $\imath(y)=(y,0)$ and $\pi(y,x)=x$, called the sequence generated by $\Phi$. 
Every short exact sequence of Banach spaces arises in this way, up to equivalence in $\Ext(X,Y)$, and the equivalence can be described in quasilinear terms as: $\Phi$ and $\Psi$ generate equivalent extensions if and only if $\Phi -\Psi = B + L$ with $B: X\to Y$ bounded and  $L: X\to Y$ linear.
 We also need the construction of the $\co(X)$ space and the universal 1-linear map $\mho: X\to \co(X)$. The construction that best suits to our present purposes is presented in \cite{cm-isr} and is slightly different from (although equivalent to) the case $p=1$ of \cite[\S~3.10]{hmbst}. Let $\mathscr H$ be a Hamel basis of $X$. Then there exists a Banach space $\co(X)$ and a 1-linear map $\mho: X\to \co(X)$, with $Q_1(\mho)= 1$ such that for every Banach space $Y$ and every 1-linear map $\Phi: X\to Y$ vanishing on $\mathscr H$ there exists a unique operator $\phi:\co(X)\to Y$ such that $\Phi=\phi\circ\mho$ with $\|\phi\|=Q_1(\Phi)$. It is clear that $\co(X)$ is spanned by the ``vectors'' $\mho(x)$ with $x\in X$ (in the sense that no proper closed subspace of $\co(X)$ contains them) and that the dual of $\co(X)$ agrees with the space of 1-linear functions $\Phi:X\to \mathbb K$ that vanish on $\mathscr H$  normed  by $Q_1(\cdot)$, and duality given by $\langle \Phi, \mho(x)\rangle = \Phi(x)$.

For every quasilinear map $\Phi:X\to \K$ one obviously has
\begin{align*}
Q_1(\Phi)&=\sup_n \left\{ \Bigg{|}\Phi\left(\sum_{i\leq n} x_i\right)- \sum_{i\leq n} \Phi(x_i )\Bigg{|}: \sum_{i\leq n} \|x_i\| \leq 1 \right\}\\
&=\sup_n \left\{ \Bigg{|}\left\langle \Phi, \mho\left(\sum_{i\leq n} x_i\right)- \sum_{i\leq n} \mho (x_i )\right\rangle \Bigg{|}: \sum_{i\leq n} \|x_i\| \leq 1 \right\}.
\end{align*}
On the other hand the universal map $\mho$, being 1-linear, generates a ``universal exact sequence''
$$
\xymatrix{0 \ar[r] & \co(X)\ar[r] & \co(X) \oplus_\mho X  \ar[r] & X  \ar[r]& 0}
$$
that in many respects behaves as a projective presentation of $X$; see below.

\medskip


We now present an intrinsic characterization of compact extensions which is almost obvious, once one does something to make it obvious. Let us agree on the following notations: given a homogeneous map $\Psi:A\to B$ acting between Banach spaces we set
$$
\nabla_\Psi=\bigcup_{n\geq 2}\left\{ \Psi\left(\sum_{i\leq n} x_i\right)- \sum_{i\leq n} \Psi(x_i ): \sum_{i\leq n} \|x_i\| \leq 1 \right\}.
$$
It is clear that $\Psi$ is 1-linear if and only if $\nabla_\Psi$ is bounded.
We say that $\Psi$ is $\nabla$-compact if $\nabla_\Psi$ is
a relatively compact subset of $B$. Note that  homogeneous $\nabla$-compact maps are automatically 1-linear: compact sets are bounded, aren't they?

\begin{prop}\label{prop:compext}
Let $\mathscr Z$ be an exact sequence $0\to Y\to Z\to X\to 0$ of Banach spaces. The following are equivalent:
\begin{itemize}
\item[(a)] $\mathscr Z$ is compact.
\item[(b)] The quotient map admits a $\nabla$-compact homogeneous section $\sigma:X\to Z$.
\item[(c)] $\mathscr Z$ is equivalent to the sequence generated by a $\nabla$-compact  map $\Psi: X\to Y$.
\item[(d)]
There is a $1$-linear map $\Phi: X\to Y$ that vanishes on $\mathscr H$ and generates a sequence equivalent to $\mathscr Z$, whose associated operator $\phi: \co(X)\to Y$ is compact.
\end{itemize}
\end{prop}
\begin{proof} Let us first prove (a)$\implies$(b). If $\mathscr Z$ is equivalent to a compact extension then there is a commutative diagram with exact rows
\begin{equation}
\xymatrix{
0 \ar[r] & K  \ar[r] \ar[d]_\alpha & P \ar[r] \ar[d]^a & X \ar[r] \ar@{=}[d] & 0 \\
0 \ar[r] & Y \ar[r] & Z \ar[r]& X \ar[r]& 0}
\end{equation}
where the upper row is a projective presentation of $X$ and $\alpha\in\fK$. If $\varsigma$ is any bounded homogeneous section of $P\to X$, then the required section of $Z\to X$ is $\sigma = a \varsigma$. The implication (b)$\implies$(c) is trivial because we can obtain the 1-linear map in (c) as $\Psi=\sigma-L$, where $L:X\to Z$ is a linear (in general unbounded) section of $Z\to X$ in which case $\nabla_\Psi$ (or, more precisely, its image in $Z$) agrees with $\nabla_\sigma$.

Assume that (c) holds and let $\Lambda:X\to Y$ be the only linear map that agrees with $\Psi$ on $\mathscr H$. Then $\Phi=\Psi-\Lambda$ is a compact 1-linear map that generates a sequence equivalent to $\mathscr Z$. Let  $\phi:\co(X)\to Y$ be  the only operator such that $\Phi=\phi\circ\mho$.
Then
 $\phi$ is compact because it maps $\nabla_\mho$ to $\nabla_\Phi=\nabla_\Psi$, a relatively compact subset of $Y$ and
therefore it maps the unit ball of $\co(X)$ (which agrees with the closed convex hull of $\nabla_\mho$,  as a
straightforward application of the Hahn--Banach separation theorem reveals) to the closed convex hull of
hull of $\nabla_\Phi$, which is compact according to a celebrated result of Grothendieck.

 Finally, to prove (d)$\implies$(a) it suffices to check that the sequence generated by $\Phi$ is itself compact. A quick look at the diagram
$$
\xymatrix{0 \ar[r] & \co(X)\ar[d]_{\phi} \ar[r] & \co(X) \oplus_\mho X  \ar[r]\ar[d]_{\text{induced}} & X  \ar@{=}[d] \ar[r]& 0\\
0 \ar[r] & Y  \ar[r]  & Y\oplus_\Phi X \ar[r]  & X \ar[r]  &0}
$$
shows that the lower row is a compact extension, which is (a).\end{proof}

The interplay between quasilinear and $1$-linear maps becomes interesting at this point. Recall from \cite[Definition 3.4.1]{hmbst} that a Banach space $X$ is said to be a $\mathscr K$-space if every quasilinear map $X\to \mathbb K$ is $1$-linear. One has:

\begin{cor} If $X$ is a $\mathscr K$-space then one can replace $n\geq 2$ by $n=2$ the definition of $\nabla_\Psi$.\end{cor}
\begin{proof} The fact that every quasilinear map $\Phi: X\to\K$ is 1-linear forces the existence a constant $M$ such that for every homogenenous $\Phi: X\to\K$ one has
$$
Q_1(\Phi)\leq M \sup\{ |\Phi(x_1+x_2)-\Phi(x_1)-\Phi(x_2)|: \|x_1\|+\|x_2\|\leq 1\}.
$$
Thus, a peaceful assimilation of the proof above implies to realize that when $X$ is a $\mathscr K$-space ``with constant $M$'' the closed convex hull of the set
$$
\{\mho(x_1+x_2)-\mho(x_1)-\mho(x_2)|: \|x_1\|+\|x_2\|\leq 1\}
$$
contains a neighborhood of the origin in $\co(X)$ --- actually the ball radius $1/M$. Probably nothing more needs to be said.\end{proof}

Let $\Phi:X\to Y$ be a 1-linear map. Following \cite[Section 3.6]{hmbst} we define
$$
Q_1[\Phi]= \inf_{\Psi\sim \Phi} Q_1(\Psi) =
 \inf_{B} Q_1(\Phi+ B),
$$
where $B$ runs over all homogeneous bounded maps from $X$ to $Y$. The function $\Phi\mapsto Q_1[\Phi]$ defines a seminorm and so a topology for 1-linear maps that is natural in many respects. 
The reader is referred to \cite[Section 4.5]{hmbst} for a variety of examples and counterexamples concerning $Q_1[\cdot]$, in particular that $
Q_1[\Phi]=0$ does not imply the triviality of $\Phi$ (Proposition 4.5.5 there). The following remark actually provides a large family of counterexamples:

\begin{prop}
Let $\Phi:X\to Y$ be a $\nabla$-compact map. If $Y$ has the BAP, then $Q_1[\Phi]=0$.
\end{prop}	

\begin{proof}[Sketch]
We have a commutative diagram (ignore the dotted arrow which represents a quasilinear map)
$$
\xymatrix{0 \ar[r] & K  \ar[r] \ar[d]_\alpha & P \ar[r] \ar[d] & X  \ar@{..>}@/_1pc/[ll]_(.3){\Omega}\ar[r] \ar@{=}[d] & 0 \\
0 \ar[r] & Y  \ar[r]  & Y\oplus_\Phi Z \ar[r] & X \ar[r] & 0}
$$
in which the upper row is a projective presentation of $X$ and $\alpha\in\fK$. The map $\Omega:X\to K$ is any quasilinear (actually 1-linear) map associated to the upper sequence and can be obtained as the difference of two sections of the quotient map $P\to X$, one bounded and the other linear. 
Since $\Phi\sim \alpha\circ\Omega$, it suffices to see that $Q_1[\alpha\circ\Omega]=0$. Pick $\varepsilon>0$ and then a finite rank $f:K\to Y$ such that $\|\alpha-f\|<\varepsilon$. Then $f\circ\Omega$ is trivial and since
$\alpha\circ\Omega= (\alpha-f)\circ\Omega+ f\circ\Omega$ one has
$$Q_1[\alpha\circ\Omega]\leq Q_1((\alpha-f)\circ\Omega)\leq
\|\alpha-f\|Q_1(\Omega)\leq \varepsilon Q_1(\Omega).\qedhere$$
\end{proof}

We do not know whether the hypothesis on the BAP can be dropped. This result, in tandem with Proposition~\ref{prop:nsingular}, reinforces the rather vague idea that compact extensions tend to be ``asymptotically trivial'', whatever that means.

None of the authors take seriously the possibility that the converse of the preceding proposition holds. Possible counterexamples are the Kalton--Peck maps $\Phi_\alpha :\ell_p\to\ell_p$ based on the Lipschitz functions $\theta_\alpha: \R^+\to\R$ defined by $\theta_\alpha(t)=\min(t,t^\alpha)$, where $1<p<\infty$ and $0<\alpha<1$. These maps generate  ``strictly nonsingular'' sequences, with the meaning that
every infinite-dimensional subspace of the quotient space
contains a further infinite-dimensional subspace on which the
quotient map is invertible. So these sequences might be counterexamples for the converse of Proposition~\ref{prop:nsingular}.

\subsection{A natural equivalence of $\Ext_{\OP K}(X,Y)$}
The same ideas can be used to obtain a nonlinear description of $\Ext_{\fK}$. We need the following observation:
\begin{lem}
If the sequences $ 0\to K \to Z \to X\to 0$ and $ 0\to K' \to Z' \to X\to 0$ are semiequivalent, then one has isomorphisms
$$
\frac{\fK(K, Y)}{\fK\text{\rm -extensible}} =
\frac{\fK(K', Y)}{\fK\text{\rm -extensible}}.
$$
for all Banach spaces $Y$.
\end{lem}

\begin{proof}
The hypothesis means that there exist a commutative diagram
$$\xymatrix{K \ar[d]_{a} \ar[r] & Z \ar[r] \ar[d] & X \ar@{=}[d]\\
K' \ar[d]_{a'}\ar[r]& Z' \ar[r]\ar[d]& X \ar@{=}[d]\\
			 K\ar[r]& Z \ar[r] & X }
$$
with short exact rows, which in turn implies that $a'a-{\bf I}_K$ is extensible to $Z$.

The composition operator $(a')^*:\mathfrak{L}(K,Y)\to \mathfrak{L}(K,Y)$ maps compact operators to compact operators and $\fK$-extensible operators to $\fK$-extensible operators and so it induces a linear map
$$
(a')^*: \frac{\fK(K, Y)}{\fK\text{\rm -extensible}} \To
\frac{\fK(K', Y)}{\fK\text{\rm -extensible}}.
$$
One also has a similar map
$$
a^*: \frac{\fK(K', Y)}{\fK\text{\rm -extensible}} \To
\frac{\fK(K, Y)}{\fK\text{\rm -extensible}}.
$$
Using the extensible character of $a'a-{\bf I}_K$ one easily realizes that $a^*(a')^*= {\bf I}_K^*$ is the identity on
$\fK(K, Y)/\fK\text{\rm -extensible}$ and reversing the roles of $a$ and $a'$ one obtains that
$(a')^*a^*$ is the identity on
$\fK(K', Y)/\fK\text{\rm -extensible}$.
\end{proof}

The universal sequence is semiequivalent to any projective presentation of $X$ and
can healthly replace it to compute $\Ext_{\fK}(X, Y)$ because, according to the lemma, one has natural isomorphisms
$$
\frac{\fK(\co X, Y)}{\fK\text{-extensible}} =
\frac{\fK(K, Y)}{\fK\text{-extensible}}.
$$
In the following result we say that a homogenenous map $\Phi:X\to Y$ is
  $\OP K$-trivial if $\Phi = K + L$ with $L$ linear and $K$ a compact map (meaning that it takes the unit ball of $X$ to a relatively compact subset of $Y$).

\begin{prop}
For all Banach spaces $X, Y$ one has isomorphisms
\begin{equation}
\Ext_{\fK}(X,Y)=\frac{\nabla\text{\rm -compact maps }X\to Y}{\fK\text{\rm -trivial maps}}
\end{equation}
\end{prop}

\begin{proof}
Note that for every map $\Phi: X\to Y$ there is a unique linear map $L: X\to Y$ such that $\Phi+L$ vanishes on a given Hamel basis of $X$. Moreover, adding a linear map to $\Phi$ does not change its $\nabla$-compact or $\OP K$-trivial character.
Hence, it suffices to see that
\begin{equation}
\frac{\fK(\co X, Y)}{\fK\text{-extensible}} =\frac{\nabla\text{\rm -compact maps }X\to Y \text{ vanishing on }\mathscr{H}}{\fK\text{\rm -trivial maps}}.
\end{equation}
We know from the proof of Proposition~\ref{prop:compext} that the numerators of these ratios coincide in the sense that compact  operators $\varphi: \co(X)\to Y$ correspond to  $\nabla$-compact maps $\Phi: X\to Y$  vanishing on $\mathscr{H}$.
 The proof will be complete after showing that $\varphi$ has a compact extension to $U$ if and only if $\Phi$ is $\fK$-trivial.

Note that if $\tilde\varphi: \co(X)\oplus_\mho X\to Y$ is any \emph{linear} extension of $\varphi$, then $\tilde\varphi(y,x)=\varphi(y)+ L(x)$, for some linear map $L:X\to Y$. If, moreover, $\tilde\varphi$ is a compact operator, then the map $K(x)=\tilde\varphi(\mho(x), x)$ is compact from $X$ to $Y$, since $\|(\mho(x), x)\|=\|x\|$. But,
$$
\Phi(x)= \varphi(\mho(x)) = \tilde\varphi(\mho(x), 0)
= \tilde\varphi(\mho(x), x)- \tilde\varphi(0, x) = K(x)-L(x),
$$
so that $\Phi$ is $\fK$-trivial. Conversely, if $\Phi$ is $\fK$-trivial, with decomposition $\Phi=K+L$, then
$\tilde\varphi(y,x)=\varphi(y)- L(x)$ defines a compact extension of $\varphi$.
\end{proof}

\section{Some remarks and open problems}\label{sec:misc}

There are a number of problems related to the research presented here that we have not been able to solve:


\subsection*{Injectivity and surjectivity of $\Ext_{\mathfrak K}\to \Ext$}

\begin{prob}
\emph{Find Banach spaces for which the natural map $\Ext_{\OP K}(X,Y)\to \Ext (X, Y)$ is nonzero and injective.}
\end{prob}

The problem asks for a closed subspace $K\subset\ell_1$ and a Banach space $Y$ such that not every compact operator $K\to Y$ admits an extension to $\ell_1$ but, when it does, it has a compact extension. The case $X=Y=\ell_2$ seems, as always, particularly interesting.

The material around Example~\ref{ex:posno} can be organized to provide explicit counterexamples very close to Hilbert spaces.
Start with a twisted Hilbert space, that is, a nontrivial self-extension $0\to\ell_2\to Z\to \ell_2\to0$ (see \cite{hmbst}, Sections 3.2 and 10.9 for the main examples). It is easily seen (cf. the proof of Theorem~\ref{Ext-K}) that there exist compact endomorphisms of $\ell_2$ that do not have extension to $Z$. If $\alpha$ is one of these, we may form the pushout
\begin{equation*}
\xymatrix{
0 \ar[r] & \ell_2  \ar[r] \ar[d]_\alpha & Z \ar[r] \ar[d] & \ell_2 \ar[r] \ar@{=}[d] & 0 \\
0 \ar[r] & \ell_2 \ar[r]^\imath & \PO \ar[r]& \ell_2 \ar[r]& 0}
\end{equation*}
to conclude that $\imath\alpha$ is a nonzero element in the kernel of $\Ext_{\OP K}(\ell_2,\PO)\to \Ext (\ell_2,\PO)$.
The point is that $\PO$ is also a twisted Hilbert space (``less twisted'' than $Z$ in fact) and so its is very close from being a straight Hilbert space. Actually, if $Z$ arises from a centralizer defined on a Banach sequence space, one can take $\alpha$ diagonal; see \cite[\S 3]{ideals}.

\begin{prob}
\emph{Find Banach spaces for which the natural map $\Ext_{\OP K}(X,Y)\to \Ext(X, Y)$ is nonzero and surjective.}
\end{prob}

The problem can be reformulated as: Does every operator $u:K\to Y$ admit a compact perturbation $u+\kappa$ extensible to $\ell_1$? While an affirmative answer could come from an exotic couple (in fact it cannot be ruled out that $\Ext(X, Y)$ has non-zero finite dimension, see \cite[Note~4.6.4]{hmbst}, but note that ``our'' $\Ext$ is denoted by $\Ext_{\bf B}$ there), we find the case in which $X=\ell_1, Y=c_0$ particularly interesting.

In \cite[Proposition 5.2.20]{hmbst}, the reader can find the smoking gun attesting that $\Ext(c_0, \ell_1)$ is not zero and, according to Theorem~\ref{Ext-K}, the same is true for $\Ext_{\fK}(c_0, \ell_1)$. It is not however easy to exhibit explicit examples, so we will give some indications here. We start again with a twisted Hilbert space, and embedding $e:\ell_2\to L_1$ and a quotient map $q:L_\infty\to \ell_2$, which  could well be the adjoint of $e$. Form the pushout and then pullback diagram (ignore the dotted arrows for the time being)
$$\xymatrix{\ell_2 \ar[d]_e \ar[r] & Z \ar[r] \ar[d] & \ell_2 \ar@{=}[d] \\
L_1 \ar@{=}[d] \ar[r]& \PO \ar[r]& \ell_2 \\
			 L_1\ar[r] \ar@{=}[d]& \PB_1 \ar[r]\ar[u] & L_\infty \ar[u]_q
			 \\
 L_1\ar@{..>}[r]& \PB_2 \ar[r] \ar@{..>}[u] \ar@{..>}[u] & c_0 \ar@{..>}[u]_f			
			 }
$$
The reader can find in \cite[Proof of Proposition 5.2.20]{hmbst} the proof that if the first row does not split, neither does the third row.
 Since $L_\infty$ is a $C(K)$-space, there is, according to \cite{ccs}, an embedding $f: c_0\to \ell_\infty$ for which the fourth row is not trivial. This is a compact sequence since the composition
 $qf$, as any operator $c_0\to \ell_2$, is compact.
The adjoint sequence
$$
\xymatrixrowsep{1pc}
\xymatrix{
0\ar[r] & c_0^* \ar[r]  \ar@{=}[d] &\PB_2^* \ar[r] & L_1^* \ar@{=}[d] \ar[r] & 0\\
&\ell_1&&L_\infty
}
$$
is again compact (see next section) and
does not split because $L_1$ is an ultrasummand.
Applying the ``fourth row argument'' as before one gets an extension of $c_0$ by $\ell_1$. The tricky point now is that all the elements of $\Ext(c_0, \ell_1)$ we ``explicitly'' know seem to appear this way. Hence, the question rather seems to decide whether non-compact extensions also exist. Proposition~\ref{prop:nsingular} once again brings us to the following question.


\begin{prob}\emph{Do singular
sequences $0\to	\ell_1\to Z\to c_0\to 0$ exist?}
 \end{prob}

This question has now the added interest that a singular sequence in $\Ext(c_0, \ell_1)$ will certainly be a noncompact sequence. Since we can solemnly declare that we are up to no good with this approach, curious readers can peruse a fragment of our Marauder's map in \cite[\S~5.2]{cs}.

\subsection*{Compactness and duality}

Schauder's theorem states that the adjoint of a compact operator between Banach spaces is again compact. On the other hand, the Hahn--Banach extension theorem implies that
the adjoint of an exact sequence of Banach spaces is again exact and it is clear that the dual of a projective space (direct summand of $\ell_1(\Gamma)$) is injective (direct summand of $\ell_\infty(\Gamma)$).

It follows that if
we have an element of $\Ext_{\fK}(X,Y)$, pictorially represented as
$$
\xymatrix{
0 \ar[r] & K  \ar[r]\ar[d]_\tau & P \ar[r] & X \ar[r] & 0\\
&Y\\
}
$$
with $\tau$ compact, then taking the adjoint of $\tau$ we obtain
and element of $\Ext_{\fK}(Y^*,X^*)$ pictorically represented by
$$
\xymatrix{
0 \ar[r] & X^*  \ar[r] & P^* \ar[r] & K^* \ar[r] & 0\\
&&&Y^*\ar[u]_{\tau^*}\\
}
$$
Note that if $\tau$ admits a compact extension to $P$, then $\tau^*$ has a compact lifting. This yields a well-defined natural ``duality" map $\Ext_{\fK}(X,Y) \to \Ext_{\fK}(Y^*,X^*)$ whose nature is quite mysterious:

\begin{prob}
\emph{Is $\Ext_{\fK}(X,Y) \to \Ext_{\fK}(Y^*,X^*)$ injective, surjective or an isomorphism for all Banach spaces $X, Y$?}
\end{prob}
	
In general, the standard duality map $\Ext(X,Y) \to \Ext(Y^*,X^*)$ is neither injective (nontrivial sequences can have trivial dual sequences) nor surjective (there are sequences formed by duals that are not duals of any sequence, see \cite{hmbst}). This problem is motivated by the fact that if $\alpha:A\to B$ is a compact operator then the range of its double adjoint
$\alpha^{**}:A^{**}\to B^{**}$ stays in $B$, which somehow suggests that the natural map $\Ext_{\fK}(X,Y) \to \Ext_{\fK}(X^{**}, Y^{**} )$ tends to be injective.



\end{document}